\documentclass[10pt,twocolumn,twoside]{IEEEtran} 
\usepackage{lscape}
\usepackage{epsfig,color,amsmath,cite}
\usepackage{amsthm}
\allowdisplaybreaks[4]
\usepackage{balance}

\IEEEoverridecommandlockouts
\usepackage{wrapfig}
\usepackage{bm}
\usepackage{epstopdf}
\usepackage{amssymb}
\usepackage{url}
\usepackage{enumitem}
\usepackage{multirow}
\usepackage{booktabs}
\usepackage{algorithm,algorithmic}	
\setlist[itemize]{leftmargin=*}

\DeclareMathOperator{\diag}{diag}

\DeclareMathOperator*{\minimize}{minimize}

\DeclareMathOperator{\subjectto}{subject\ to}
\DeclareMathOperator{\vecc}{vec}
\makeatother

\newtheorem{mydef}{Definition}

\newtheorem{rem}{Remark}
\newtheorem{asmp}{Assumption}
\newtheorem{prop}{Proposition}
\newtheorem{corol}{Corollary}
\newtheorem{lem}{Lemma}

\newcommand{\mat}[1]{\boldsymbol{#1}}

\newcommand{\normof}[1]{\|#1\|}

\newcommand{\bmat}[1]{\begin{bmatrix} #1 \end{bmatrix}}

\providecommand{\mA}{\ensuremath{\mat{A}}}
\providecommand{\mB}{\ensuremath{\mat{B}}}
\providecommand{\mC}{\ensuremath{\mat{C}}}
\providecommand{\mD}{\ensuremath{\mat{D}}}
\providecommand{\mE}{\ensuremath{\mat{E}}}

\providecommand{\mG}{\ensuremath{\mat{G}}}
\providecommand{\mH}{\ensuremath{\mat{H}}}
\providecommand{\mI}{\ensuremath{\mat{I}}}

\providecommand{\mK}{\ensuremath{\mat{K}}}
\providecommand{\mL}{\ensuremath{\mat{L}}}

\providecommand{\mO}{\ensuremath{\mat{O}}}

\providecommand{\mQ}{\ensuremath{\mat{Q}}}
\providecommand{\mR}{\ensuremath{\mat{R}}}
\providecommand{\mS}{\ensuremath{\mat{S}}}

\providecommand{\mV}{\ensuremath{\mat{V}}}

\providecommand{\mY}{\ensuremath{\mat{Y}}}
\providecommand{\mZ}{\ensuremath{\mat{Z}}}

\newcommand{\st}{{\rm s.t.}}

\usepackage{refcount}

\newcommand{\myfootnotetext}[1]{\footnotetext{#1\label{fn:text}%
		\edef\fnmark{\getpagerefnumber{fn:mark}}%
		\edef\fntext{\getpagerefnumber{fn:text}}%
		\ifx\fnmark\fntext\else\ClassWarning{}{footnote mark and text on different pages!}\fi}}

\newcommand{\m}{\boldsymbol}
\allowdisplaybreaks[4]
\usepackage[colorlinks = true,
linkcolor = blue,
urlcolor  = blue,
citecolor = blue,
anchorcolor = blue]{hyperref}

\usepackage{tikz}
\usetikzlibrary{positioning,fit,shapes.geometric,backgrounds}

\usetikzlibrary{arrows.meta}
\tikzset{%
	>={Latex[width=1mm,length=1mm]},
	base/.style = {rectangle, rounded corners, draw=black,
		minimum width=2cm, minimum height=0.5cm,
		text centered},
	activityStarts/.style = {base, fill=white!30},
	startstop/.style = {base, fill=white!30},
	activityRuns/.style = {base, fill=white!30},
	process/.style = {base, fill=white!30},
}

\tikzstyle{block} = [draw, fill=white!20, rectangle, 
minimum height=6em, minimum width=6em]
\tikzstyle{sum} = [draw, fill=blue!20, circle, node distance=1cm]
\tikzstyle{input} = [coordinate]
\tikzstyle{output} = [coordinate]
\tikzstyle{pinstyle} = [pin edge={to-,thin,black}]

\tikzset{%
	>={Latex[width=1mm,length=1mm]},
	base/.style = {rectangle, rounded corners, draw=black,
		minimum width=2cm, minimum height=0.5cm,
		text centered},
	activityStarts/.style = {base, fill=white!30},
	startstop/.style = {base, fill=white!30},
	activityRuns/.style = {base, fill=white!30},
	process/.style = {base, fill=white!30},
}

\usetikzlibrary{shapes,arrows}
\usetikzlibrary{calc}

\title{Time-Varying Sensor and Actuator Selection for Uncertain Cyber-Physical Systems}
\author{Ahmad F. Taha, \IEEEmembership{Member,~IEEE,} Nikolaos Gatsis~\IEEEmembership{Member,~IEEE,}  Tyler Summers~\IEEEmembership{Member,~IEEE}, Sebastian Nugroho,~\IEEEmembership{Student Member,~IEEE. }\\
	\thanks{
Ahmad F. Taha, Nikolaos Gatsis, and Sebastian Nugroho are affiliated with the Department of Electrical and Computer Engineering at the University of Texas San Antonio. Tyler Summers is with the Department of Mechanical Engineering, University of Texas at Dallas. Emails: \{ahmad.taha, nikolaos.gatsis,sebastian.nugroho\}@utsa.edu, tyler.summers@utdallas.edu. This material is based upon work supported in part by the National Science Foundation under Grants No. ECCS-1462404, CMMI 1728629, and CMMI 1728605. The work of T. Summers was partially sponsored by the Army Research Office and was partially accomplished under Grant Number: W911NF-17-1-0058. }}
\begin{document}
	\maketitle
	\begin{abstract}
We propose methods to solve time-varying, sensor and actuator (SaA) selection problems for uncertain cyber-physical systems. We show that many SaA selection problems for optimizing a variety of control and estimation metrics can be posed as semidefinite optimization problems with mixed-integer bilinear matrix inequalities (MIBMIs). {Although this class of optimization problems are computationally challenging, we present tractable approaches that directly tackle MIBMIs, providing both upper and lower bounds, and that lead to effective heuristics for SaA selection.} The upper and lower bounds are obtained via successive convex approximations and semidefinite programming relaxations, respectively, and selections are obtained with a novel slicing algorithm from the solutions of the bounding problems. Custom branch-and-bound and combinatorial greedy approaches are also developed for a broad class of systems for comparison. Finally, comprehensive numerical experiments are performed to compare the different methods and illustrate their effectiveness. 
\end{abstract}

\begin{IEEEkeywords}
Sensor and actuator selection, cyber-physical systems, linear matrix inequalities, controller design, observer design, mixed integer programming.
\end{IEEEkeywords}

%\vspace*{-0.5cm}
\section{Introduction \& Brief Literature Review}
\IEEEPARstart{M}{any} emerging complex dynamical networks, from critical infrastructure to industrial cyber-physical systems (CPS) to various biological networks, are increasingly able to be instrumented with new sensing and actuation capabilities. These networks comprise growing webs of interconnected feedback loops and must operate efficiently and resiliently in dynamic and uncertain environments. The prospect of incorporating large numbers of additional sensors and actuators (SaAs) raises fundamental and important problems of jointly and dynamically selecting the most effective SaAs, in addition to simultaneously designing corresponding estimation and control laws associated with the selected SaAs.

% brief literature review, and SHARP distinction with present work
There are many different quantitative notions of network controllability and observability that can be used as a basis for selecting effective SaAs in uncertain and dynamic cyber-physical networks. Notions based on classical Kalman rank conditions for linear systems focus on binary structural properties \cite{liu2011controllability,nepusz2012controlling,ruths2014control,olshevsky2014minimal,pequito2016}. More elaborate quantitative notions based on Gramians \cite{pasqualetti2014controllability,summers2014optimal,summers2014submodularity,tzoumas2016,yannature2015,zhao2016scheduling,nozari2016time,Chanekar2017,jadbabaie2018deterministic} and classical optimal and robust control and estimation problems \cite{Polyak-LMI_sparse_fb,Dhingra2014,munz2014sensor,Argha2016,summers2016actuator,summers2017performance,zhang2017sensor,Taylor2017} for linear systems have also been studied. For selecting SaAs based on these metrics, several optimization methods are proposed in this literature, including combinatorial greedy algorithms~\cite{tzoumas2016,tzoumas2016near,zhang2017sensor,summers2014submodularity,summers2016actuator}, convex relaxation heuristics using sparsity-inducing $\ell_1$ penalty functions~\cite{Polyak-LMI_sparse_fb,Dhingra2014,munz2014sensor,Argha2016} and reformulations to mixed-integer semidefinite programming via the big-M method or McCormick's relaxation~\cite{Taylor2017,Chanekar2017,5717988}. As a departure from control-theoretic frameworks, the authors in~\cite{haber2018state} explore an optimization-based method for reconstructing the initial states of nonlinear dynamic systems given \textit{(a)} arbitrary nonlinear model, while \textit{(b)} optimally selecting a fixed number of sensors. The authors also showcase the scalability of their approach in comparison with sensor selection algorithm based on empirical observability of nonlinear systems. 

Despite the recent surge of interest in quantifying network controllability and observability and in associated SaA selection problems, a much wider set of metrics are relevant for uncertain cyber-physical systems. {The existing literature tends to focus mainly on classical metrics (e.g., involving Kalman rank~\cite{liu2011controllability},  Gramians~\cite{nozari2016time,summers2016actuator,tzoumas2016}, Linear Quadratic Regulators~\cite{summers2017performance,summers2016actuator,Chanekar2017}, and Kalman Filters~\cite{Taylor2017,zhang2017sensor}) and deterministic linear time-invariant systems}. Methods for time-varying systems with various uncertainties and constraints are also important to broaden applicability. It is well known that a broad variety of systems and control problems can be cast in the form of semidefinite programs (SDP) and linear matrix inequalities (LMI) \cite{boyd1994linear}, but many of these more recent formulations have not been considered in the context of SaA selection. {In general, the selection of sensors or actuators  and design of associated estimation and control laws for many metrics can be posed as semidefinite optimization problems with mixed-integer bilinear matrix inequalities (MIBMIs) as we have recently shown in~\cite{Taha2017c}. A general MIBMI formulation for the selection problem is also discussed in the ensuing sections. }

% contributions
Here we propose methods to solve time-varying, sensor and actuator (SaA) selection problems for uncertain cyber-physical systems. Our methods can be applied to any of the broad range of problems formulated as MIBMIs. Although this class of optimization problems is computationally challenging, we present tractable approaches that provide upper and lower bounds and lead to effective heuristics for SaA selection. The upper and lower bounds are obtained via successive convex approximations and SDP relaxations, respectively, and selections are obtained with a novel slicing algorithm from the solutions of the bounding problems. 
%The direct derivation of approximations and relaxations from the MIBMIs contrasts our approach with other existing optimization-based approaches. The time-varying aspect allows SaA selections to change dynamically based on time-varying conditions in the network. Mixed-integer SDPs and combinatorial greedy approaches are also developed for a broad class of systems for comparison. Finally, comprehensive numerical experiments are performed to compare the different methods and illustrate their effectiveness. 

A preliminary version of this work appeared in \cite{Taha2017c} where we developed customized algorithms for actuator selection.  Here we significantly extended the methodology with the successive convex approximation and convex relaxation approaches and provide comprehensive numerical experiments. 
%The extended version of this paper can be found in~\cite{Taha2017d}. The extended version includes significant additions to the paper: (1) Extensions of this work to include a variety of other control and estimation metrics with sensor/actuator selection. This extension is showcased in~\cite[Appendix E]{Taha2017d}. (2) An alternate formulation for the SaA selection problem via the big-M method that amounts to solving mixed-integer SDPs, in comparison with the convex relaxations and approximations we develop in this paper. (3) A thorough discussion on greedy algorithms, and extended numerical tests for another dynamical system with different network size and properties. 

 \textit{Paper Notation} --- Italicized, boldface upper and lower case characters represent matrices
and column vectors: $a$ is a scalar, $\m a$ is a vector, and $\m A$ is a matrix.  Matrix $\m I_n$ is the identity square matrix of size $n$, and vector $\m 1_n$ is a vector of ones of size $n$;  $\m O_{m \times n}$ defines a zero matrix of size $m \times n$. $\mathbb{S}^n$ denotes the set of symmetric matrices of size $n$; $\mathbb{S}^n_+$ and $\mathbb{S}^n_{++}$ are the sets of symmetric positive semidefinite and positive definite matrices. $\| \m A\|_*$ denote the nuclear norms of $\m A$. The symbol $\diag(\m a)$  denotes a diagonal matrix whose diagonal entries are given by the vector $\m a$; $\mathrm{diag}(\m A)$ forms a column vector by extracting the diagonal entries of $\m A$. The symbol $\Lambda(\m A)$ denotes the set of complex eigenvalues of a matrix $\m A$. 

The next section presents the framework and problem formulation, and details the paper contributions and organization. 
%The notation used in the paper is included in Appendix~\ref{app:notation}. 
\section{CPS Model and Paper Contributions}~\label{sec:PbmFormulation}
 We consider time-varying CPSs with $N$ nodes modeled as
 \vspace{-0.3cm}
		\begin{subequations}~\label{equ:CPSModelSaA}
		\begin{align}
		\dot{\m x}(t)  &=\mA^j\m x (t) +\mB_u^j\m \Pi^j \m u (t)  + \mB_w^j \m w (t) + \mB_f^j\m f^j(\m x),\\
		\m y (t)  &=  \m \Gamma^j \mC^j\m x(t)  + \mD_u^j \m u (t)  + \mD_v^j \m v(t), \;\; \m x^j(t_0)=\m x_0^j
		\end{align}
	\end{subequations}
The network state $\m x(t) \in \mathbb{R}^{n_x}$ consists of each of $N$ nodal agent states $\m x_i \in \mathbb{R}^{n_{x_i}}$, $i=1,\ldots,N$. Each nodal agent has a set of available inputs $\m u_i \in \mathbb{R}^{n_{u_i}}$ and measurements $\m y_i(t) \in \mathbb{R}^{n_{y_i}}$. The mapping from the input to state vector can thus be written in the form $\mB_u= \mathrm{blkdiag}(\mB_{u_1}, \ldots, \mB_{u_N}) $. The system nonlinearity can be expressed as $\m f(\m x) \in \mathbb{R}^{n_x}$ and $\m B_f$ represents the distribution of the nonlinearities. The vectors $\m w(t) \in \mathbb{R}^{n_w}$ and $\m v(t) \in \mathbb{R}^{n_v}$ model unknown inputs and data perturbations. In summary, the system has $n_x$ states, $n_u$ control inputs, $n_y$ output measurements, $n_w$ unknown inputs, and $n_v$ data perturbations. that are common in CPSs. Superscript $j$ denotes the time-period (see Remark~\ref{rem:TP}).
The model~\eqref{equ:CPSModelSaA} includes binary variables  $\pi_i$, $i=1,\ldots,N$, where $\pi_i=1$ if the actuator of the $i$-th nodal agent is selected, and 0 otherwise. Similarly, we define binary variables $\gamma_i$, $i=1,\ldots,N$, where $\gamma_i=1$ if the sensor of the $i$-th nodal agent is selected, and 0 otherwise. Variables $\pi_i$ and $\gamma_i$ are organized in vectors $\m\pi$ and $\m\gamma$, i.e., $\m\Pi=\mathrm{blkdiag}(\pi_1 \mI_{n_{u_1}}, \ldots, \pi_N \mI_{n_{u_N}})$ and $\m\Gamma=\mathrm{blkdiag}(\gamma_1 \mI_{n_{y_1}}, \ldots, \gamma_N \mI_{n_{y_N}})$. 
 
\begin{rem}[Topological Evolution]~\label{rem:TP}
	In~\eqref{equ:CPSModelSaA}, the optimal SaA selection and the control/estimation laws change from one time-period to another. The time-frame depends on the application under study, and the state-space matrices are obtained through an apriori analysis of the system dynamics. For example, in power networks the state-space matrices $(\m A^j, \m B_u^j, \ldots)$ change according to the operating point of the system which is determined via optimal power flow routines~\cite{dommel1968optimal}. The time-horizon of this change is around 5 minutes. In water distribution networks, this change is often in hours as the water demand patterns and water flows evolve at a slower time-scale than electric power demand~\cite{eliades2011fault}. In this paper, we assume that the transition in the state-space matrices is given. 
\end{rem}
 
%The high-level problem for time-varying SaA selection can be formulated as follows:
%\begin{subequations}\label{equ:HLPbm}
%	\begin{align}
%	\min\hspace{1.3cm} &   \hspace{-1.3cm}\sum_{j=1}^{T_f} \left\lbrace c_{\pi,j}(\m \pi^j)+c_{\gamma,j}(\m \gamma^j)+  \mathrm{CtrlObj} + \mathrm{EstObj}\right\rbrace \\
%	\subjectto \;\;\;&   \eqref{equ:CPSModelSaA},\;\;\{\m\pi^j\}_{j=1}^{T_f}\in \mathcal{P}, \{\m \gamma^j\}_{j=1}^{T_f} \in \mathcal{G}\\
%	& \hspace{-0.7cm} \mathrm{CtrlConstraints}(\m \pi^j), \mathrm{EstConstraints}(\m \gamma^j).
%	\end{align}
%\end{subequations}
%In~\eqref{equ:HLPbm}, the objective is to minimize the sum of the cost as a function of the selected SaAs and certain control and estimation metrics, which are described in the next section. The constraints include the dynamical system evolution for all time-periods $j \in \{1,\ldots,T_f\}$, and operational set constrains ($\mathcal{P}$ and $\mathcal{G}$) on $\m \gamma$ and $\m \pi$ that include the binary nature of these variables and possibly \textit{ramp constraints} as a function of previous SaA selection. The control and estimation constraints depend on the considered estimation and control metrics which can be written as SDPs in many cases as discussed in the next section. 

%\vspace{-0.4cm}

%\subsection{Paper Organization  and Contributions}
{The formulations in this paper are building on semidefinite programming (SDP) approaches for robust control and estimation routines; see ~\cite{boyd1994linear,vanantwerp2000tutorial}.
To set the stage, we first succinctly list control and estimation formulations as SDPs in Appendix~\ref{app:cont}, where the system dynamics, controller/observer form, optimization variables, and the optimization problem are stated.  The listed formulations are instrumental in formalizing the SaA selection problem since the LMIs share a similar structure. Many other control and estimation laws can fit directly into the proposed methodologies. The main  contributions of this paper are detailed next.}
\begin{itemize}
\item   First, we show that a large array of optimal control and estimation problems with SaA selection share a similar level of computational complexity of solving optimization problems with MIBMIs (Section~\ref{sec:BeyondStab}).
\item Second, we develop one-shot convex relaxation that produces a lower bound to the original problem with MIBMIs. With respect to previous SDP relaxations~\cite{BV-SDR-Control}, the proposed approach entails matrix variables of smaller dimension, which is computationally advantageous. Two successive convex approximations that yield upper bounds are also developed. Theoretical guarantees on the convergence of the convex relaxations and approximations are provided, with the necessary background and assumptions. The successive convex approximations draw from previous general methods~\cite{dinh2012combining, LeeHu2016}, but this paper develops specialized algorithms for the MIBMI problem structures that stem specifically from sensor and actuator selection. We also develop simple algorithms to recover the binary selection of SaAs, in addition to the state-feedback gains and performance indices  (Sections~\ref{sec:introBMIs}--\ref{sec:succapproxMIBMIs}).
\item Third, we include a general formulation that utilizes the big-M method, thereby transforming the optimization problem that includes MIBMIs to a mixed-integer semidefinite program (MISDP)---this approach is detailed Section~\ref{sec:bigM}. 
%In addition, classical greedy algorithms are explored in Section~\ref{sec:Suboptimal}.
\end{itemize}
Comprehensive numerical examples are provided in Section~\ref{sec:numtests}. The results show the performance of the developed methods and that the optimal solution to the relaxed MIBMIs is nearly obtained in mostly all instances of our study. The numerical results also corroborate the theoretical results, and the necessary assumptions needed to obtain convergence are satisfied. 
\normalcolor
The next section presents the developed framework of time-varying SaA selection for uncertain dynamic systems. 
\section{Time-Varying SaA Selection with Various Metrics: A Unifying MIBMI Framework}
	\label{sec:BeyondStab}
 	In this section, we show that a plethora of control or estimation problems with time-varying SaA can be written as nonconvex optimization problems with MIBMIs. This observation considers different formulations pertaining to various observability and controllability metrics.
In particular, replacing $\m B_u$ with $\m B_u \m \Pi$ and $\m C$ with $\m \Gamma\m C$ in the SDPs in Appendix~\ref{app:cont} significantly increases the complexity of the optimization problem. This transforms the SDPs into nonconvex problems with MIBMIs, thereby necessitating the development of advanced optimization algorithms---the major contribution of this paper. 

For concreteness, we only consider the actuator selection problem for robust $\mathcal{L}_{\infty}$ control of uncertain linear systems (see the second row of Appendix~\ref{app:cont}  or~\cite{pancake2000analysis}), and leave the other SDP formulations with different control/estimation metrics as simple extensions.
%For systems dynamics $\dot{\m x}(t) = \m A\m x(t) + \m B_u \m u(t) + \m B_w \m w(t)$ and performance index $\m z(t) = \m C \m x(t) + \m D_{wz}\m w(t)$, the $\mathcal{L}_{\infty}$ control with full set of actuators 
\normalcolor
Under this simplifying setup and focusing on the robust control with actuator selection, we can write the system dynamics as:
	\begin{subequations}~\label{equ:CPSModelUIsDAs}
	\begin{align}
	\dot{\m x}(t)  &=\mA^j\m x(t)  +\mB_u^j\m \Pi^j \m u(t)   + \mB_w^j \m w(t)  \\
	\m z(t) & = \m C_z^j \m x(t) + \m D_{wz}^j \m w(t),
	\end{align}
	\end{subequations}
	where $\m\Pi^j$ is binary matrix variable (cf. Section~\ref{sec:PbmFormulation}) and $\m z(t)$ is the control performance index. The time-varying sequence of selected actuators and stabilizing controllers is obtained as the solution of the following multi-period optimization problem:
	\begin{subequations}
				\label{eq:StabUI-all}
				\begin{align}
				\minimize_{\{\mS,\mZ, \zeta, \m\pi\}^j} & \;\;\; \sum_{j=1}^{T_f} (\eta+1)\zeta^j +   \alpha_{\pi}^{\top} \m \pi^j \\
\hspace{-0.2cm}\subjectto 				& \bmat{\mA^j\mS^j+\mS^j\mA^{j\top}+\alpha\mS^j \\-\mB_u^j \m\Pi^j \mZ^j -\mZ^{j\top} \m\Pi^j \mB_u^{j\top} &\mB_w^j \\\mB_w^{j\top} & -\alpha\eta \mI} \preceq \m {O} \label{eq:StabUI-MIBMI}\\
				& \bmat{-\mS^j & \mO & \mS^j\mC_z^{^j\top}\\
					\mO & -\mI & \mD_{wz}^{\top}\\
					\mC_z^j\mS^j& \mD_{wz}^j &-\zeta^j\mI} \preceq \m {O}\\
& 	\m H \m \pi \leq \m h,\; \m \pi \in  \{0,1\}^N.
				\end{align}
			\end{subequations}		
In~\eqref{eq:StabUI-all}, the optimization variables are matrices $(\mS, \mZ, \m Y)^j$, the actuator selection $\m \pi^j$ (collected in vector $\m \pi$ for all $j$), and the robust control index $\zeta^j$ for all $j \in \{1,\ldots,T_f\}$, where $\alpha$ and $\eta$ are predefined positive constants~\cite{pancake2000analysis}. Given the solution to \eqref{eq:StabUI-all}, the stabilizing control law for the $\mathcal{L}_{\infty}$ problem can be written as $\m u^*(t)=-\mZ^{*j}(\m S^{*j})^{-1}\m x(t)$ for all $t \in [t_{j},t_{j+1})$. This guarantees that $\normof{\m z(t)}_2 \leq \sqrt{(\eta+1)\zeta^*} \normof{\m w(t)}_{\infty}$. The constraint $\m H\m \pi \leq \m h$ couples the selected actuators across time periods,
and is a linear logistic constraint that includes the scenarios discussed in Appendix~\ref{app:log}.
%\begin{itemize} 
%	\item Activation and deactivation of SaAs in a specific time-period $j$. For example, if actuator $i$ cannot be selected at period $j$, we set $\pi_i^j \leq 0$.
%	\item If actuator $k$ is allowed to be selected only after actuator $i$ is selected at period $j$, we set
%	$\pi_k^{j+1} \leq \pi_i^j,$ for $j=1,\ldots, T_f.$
%	\item  If actuator $k$ must be deselected after actuator $i$ is selected at period $j$, we set
%	$\pi_k^{j+1} \leq 1-\pi_i^j,$ for $j=1,...T_f$.
%	\item Upper and lower bounds on the total number of active SaAs per period can be accounted for.
%	\item Other constraints such as minimal number of required active actuators in a certain region of the CPS, and unit commitment constraints that are obtained from solutions day-ahead planning problems, can be included.
%	%	\footnote{In power systems, for example, this could entail that a generator or distributed energy resource (i.e., an actuator) $k$ will become inactive for a certain periods of time. }
%\end{itemize}
The optimization problem \eqref{eq:StabUI-all} includes MIBMIs due to the term $\m B_u^j \m\Pi^j\mZ^j$.
The bilinearity together with the integrality constraints bring about the need for specialized optimization methods.  It should be emphasized  that~\eqref{eq:StabUI-all} is \emph{not} a mixed-integer convex program. Therefore, general-purpose mixed-integer convex programming solvers are not applicable.

 		Interestingly, the design of the remaining controllers and observers in~Appendix~\ref{app:cont} largely share the optimization complexity of~\eqref{eq:StabUI-all}. It can be observed that \emph{all} design problems in Appendix~\ref{app:cont} feature MIBMIs with the form $\mB_u \m\Pi \mZ +\mZ^{\top} \m\Pi \mB_u^{\top}$ or a similar one. This simple idea signifies the impact of finding a solution to optimization problems with MIBMIs. In fact, many LMI formulations for control problems in~\cite{boyd1994linear} become MIBMIs when SaA selection is included. Using~\eqref{eq:StabUI-all} as an exemplification for other problems with similar non-convexities,  custom optimization algorithms to deal with such MIBMIs are proposed in the ensuing sections.  
\vspace{-0.2cm}

\section{From MIBMIs to BMIs}
\label{sec:introBMIs}
{This section along with Sections~\ref{sec:relaxMIBMIs} and \ref{sec:succapproxMIBMIs} develops a series of methods to deal with MIBMIs that all have the same starting point: Relaxing the integer constraints to continuous intervals.} The resulting problem is still hard to solve, as it includes bilinear matrix inequalities (BMIs). {  For clarity, we consider a single-period version of the $\mathcal{L}_{\infty}$ problem with actuator selection, i.e., problem~\eqref{eq:StabUI-all} with $T_f=1$.} This section presents some preparatory material that will be useful in the next sections. 
We start by considering the actuator selection problem with optimal value denoted by $f^*$.
\begin{subequations}
				\label{eq:StabUI2-all}
				\begin{align}
				f^*= \minimize_{\mS, \mZ,  \zeta, \m\pi} & \;\;\;\;\;(\eta+1)\zeta + \m \alpha_{\pi}^{\top} \m \pi \label{eq:StabUI2-OF} \\
\subjectto 		& \bmat{\mA \mS +\mS \mA^{\top}+\alpha\mS  \\-\mB_u  \m\Pi  \mZ  -\mZ^{\top} \m\Pi  \mB_u^{\top} &\mB_w  \\\mB_w^{\top} & -\alpha\eta\mI} \preceq \m {O} \label{eq:StabUI2-MIBMI}\\
				&  \bmat{-\mS  & \mO & \mS \mC_z^{ \top}\\
					\mO & -\mI & \mD_{wz}^{\top}\\
					\mC_z \mS & \mD_{wz}  &-\zeta \mI} \preceq \m {O} \label{eq:StabUI2-zetaS}\\
					&  \m H \m\pi \leq \m h\label{eq:StabUI2-coupl}\\
					& \m\pi  \in \{0,1\}^N. \label{eq:StabUI2-int}
				\end{align}
			\end{subequations}		
			
The following standing assumption regarding the feasibility of~\eqref{eq:StabUI2-all} is made throughout the paper.
\begin{asmp}
Problem~\eqref{eq:StabUI2-all} is feasible for $\pi_i=1$, $i=1,\ldots,N$ with constraints~\eqref{eq:StabUI2-MIBMI}, \eqref{eq:StabUI2-zetaS}, and~\eqref{eq:StabUI2-coupl} holding as strict inequalities.
\label{asmp:Feas}
\end{asmp}
The previous assumption essentially postulates that when all actuators are selected, then  $\mS, \mZ,  \zeta$ can be found so that matrix inequalities~\eqref{eq:StabUI2-MIBMI} and~\eqref{eq:StabUI2-zetaS} hold with $\m {O}$ on the left-hand side replaced by $-\epsilon \m I$, and \eqref{eq:StabUI2-coupl} with $\m h$ replaced by $\m h-\epsilon'\bm{1}$, for sufficiently small $\epsilon>0$ and $\epsilon'>0$. Such a point does \emph{not} need to be the optimal solution~\eqref{eq:StabUI2-all}; Assumption~\ref{asmp:Feas} only requires the existence of such a point in the feasible set.  It follows from the previous discussion that finding such a point is a convex optimization problem. 

The methods developed in Sections~\ref{sec:relaxMIBMIs} and \ref{sec:succapproxMIBMIs} rely on substituting the integer constraint~\eqref{eq:StabUI2-int} with the box constraint
\begin{equation}
\m 0 \leq \m\pi  \leq \m 1.   \label{eq:StabUI2-box}
\end{equation}
{  Problem~\eqref{eq:StabUI2-all} with~\eqref{eq:StabUI2-int} substituted by~\eqref{eq:StabUI2-box} can be written as}
\begin{subequations}
\label{eq:NLSDP-all}
\begin{align}
L=\minimize_{\m p} ~~~& f(\m p) \label{eq:NLSDP-obj} \\
\subjectto~~~& \mathcal{G}(\m p) \preceq \m O \label{eq:NLSDP-MI}
\end{align}
\end{subequations}
{  where the shorthand notation  $\m p= [\vecc(\m S)^{\top} \,\, \zeta \,\, \vecc(\m Z)^{\top}\,\,  \bm\pi)^{\top}]^{\top}$ is used to denote the optimization variables.} The objective is $f(\m p)=\zeta  + \m \alpha_{\pi}^{\top} \m \pi$, and $\mathcal{G}(\m p)$ is a matrix-valued function that includes the left-hand sides of~\eqref{eq:StabUI2-MIBMI}, \eqref{eq:StabUI2-zetaS}, \eqref{eq:StabUI2-coupl}, and the two sides of~\eqref{eq:StabUI2-box}, in a block diagonal form. Problem~\eqref{eq:NLSDP-all} has the general form of a nonlinear SDP~\cite{shapiro1997first}. The dimensions of $\m p$ and $\mathcal{G}(\m p)$ are respectively given by $\m p\in\mathbb{R}^d$ and $\mathcal{G}(\m p)\in \mathbb{S}^{\kappa}$, where $d$ and $\kappa$ can be inferred from~\eqref{eq:StabUI2-all}. The notation $D\mathcal{G}(\m p)$ is used for the differential of $\mathcal{G}(\m p)$ at $\m p$, i.e.,  $D\mathcal{G}(\m p)$ maps a vector $\m q \in \mathbb{R}^d$ to $\mathbb{S}^\kappa$ as follows
\begin{equation}
\label{eq:Gdiff}
[D\mathcal{G}(\m p)] \m q = \sum_{i=1}^d q_i \frac{\partial \mathcal{G}(\m p)}{\partial p_i}. 
\end{equation}
The optimal value serves as an index to formally compare the various formulations to be developed in the sequel. But comparison with respect to control metrics is also important, therefore, the resulting controllers are also evaluated in terms of the system closed-loop eigenvalues in the numerical tests of Section~\ref{sec:numtests}. 
The relationship between the optimal value of~\eqref{eq:StabUI2-all} and~\eqref{eq:NLSDP-all} is formalized in the following proposition.
\begin{prop}
\label{prop:intrelax}
With $L$ denoting the optimal value of problem~\eqref{eq:NLSDP-all}, it holds that $L\leq f^*$.
\end{prop}
 
\begin{IEEEproof}[Proof of Proposition~\ref{prop:intrelax}]
The proposition holds because~\eqref{eq:StabUI2-box} represents a relaxation of~\eqref{eq:StabUI2-int}. 
\end{IEEEproof}\normalcolor

Problem~\eqref{eq:NLSDP-all} is still hard to solve, because it contains the BMI~\eqref{eq:StabUI2-MIBMI}. Since the problem is nonconvex, several algorithms seek to find a stationary point of~\eqref{eq:NLSDP-all}, instead of a globally optimal one. Before formally stating the definition of stationary point, the Lagrangian function of~\eqref{eq:NLSDP-all} is given next: 
\begin{equation}
\mathcal{L}(\m p, \bm\Lambda) = f(\m p) + \mathbf{trace}[\m \Lambda \mathcal{G}(\m p)],
\label{eq:LagrFnNLSDP}
\end{equation}
where $\m \Lambda$ is a Lagrange multiplier matrix. Stationary points of~\eqref{eq:NLSDP-all} abide by the following definition.
\begin{mydef}
\label{def:KKT}
A pair $(\m p^*, \bm\Lambda^*)$ is a KKT point of~\eqref{eq:NLSDP-all}, and $\m p^*$ is a stationary point of~\eqref{eq:NLSDP-all}, if the following hold:
1) Lagrangian optimality: $\nabla_{\m p}\mathcal{L}(\m p^*,\m \Lambda)=\bm{0}$; 2) primal feasibility: $\mathcal{G}(\m p^*)\preceq \m O$; 3) dual feasibility: $\m \Lambda^* \succeq \m O$; and 4)  complementary slackness: $\m \Lambda^* \mathcal{G}(\m p^*)=\m O$.
\end{mydef}
Conditions 1)--4) in the above definition are the KKT conditions for~\eqref{eq:NLSDP-all}. {  These become necessary conditions that locally optimal solutions of~\eqref{eq:NLSDP-all} must satisfy, when appropriate constraint qualifications hold. Constraint qualifications are properties of the feasible set of an optimization problem; in particular, they are desirable conditions that the constraints of the optimization problem must satisfy.} To make this concept concrete, we present two typical constraint qualifications next~\cite{shapiro1997first}.
\begin{mydef}
\label{def:slater}
Problem~\eqref{eq:NLSDP-all} satisfies Slater's constraint qualification if there is a point $\m p^0\in\mathbb{R}^d$ satisfying $\mathcal{G}(\m p^0)\prec \m O$.
\end{mydef}
Slater's constraint qualification guarantees zero duality gap for problems of the form~\eqref{eq:NLSDP-all}  when $f(\m p)$ and $\mathcal{G}(\m p)$ are convex. Though $\mathcal{G}(\m p)$ is not convex for the problem at hand, we will use Slater's constraint qualification for convex approximations of~\eqref{eq:NLSDP-all} in the sequel.  A constraint qualification useful for nonconvex nonlinear SDPs is given next.
\begin{mydef}
\label{def:MFCQ}
The Mangasarian-Fromovitz constraint qualification (MFCQ) holds at feasible point $\m p^0$ if there exists a vector $\m q\in\mathbb{R}^d$ such that
\begin{equation}
\mathcal{G}(\m p^0) + [D\mathcal{G}(\m p^0)] \m q \prec \m O. 
\label{eq:MFCQ}
\end{equation}
\end{mydef}
Under MFCQ, the KKT conditions become necessary for local optima of~\eqref{eq:NLSDP-all}.
\begin{lem}
\label{lem:1stOrderOpt}
Let $\m p^*$ be a locally optimal solution of~\eqref{eq:NLSDP-all}. Then under MFCQ, there exists a Lagrange multiplier matrix $\m \Lambda^*$ that together with $\m p^*$ satisfies the KKT conditions of Definition~\ref{def:KKT}. 
\end{lem}
 
\begin{IEEEproof}[Proof of Lemma~\ref{lem:1stOrderOpt}]
This result is typical in the literature of nonlinear SDPs; see \cite[Sec.~4.1.3]{shapiro2000duality}.
\end{IEEEproof}\normalcolor
The significance of Lemma~\ref{lem:1stOrderOpt} is that it characterizes the points which are local minima of~\eqref{eq:NLSDP-all}.
For future use, we mention next two refinements of the KKT conditions of Definition~\ref{def:KKT}. Specifically, the complementary slackness condition implies that $\mathbf{rank}[\mathcal{G}(\m p^*)]+ \mathbf{rank}(\m \Lambda^*)\leq \kappa$~\cite[p.~307]{shapiro1997first}. A stricter condition is defined as follows.
\begin{mydef}
\label{def:strictcompl}
A KKT point of~\eqref{eq:NLSDP-all} satisfies the strict complementarity if $\mathbf{rank}[\mathcal{G}(\m p^*)]+ \mathbf{rank}(\m \Lambda^*)= \kappa$.
\end{mydef}
To state the second condition, the definition of a feasible direction for problem~\eqref{eq:NLSDP-all} is provided next.
\begin{mydef}
Let $\m p^0$ be a feasible point of~\eqref{eq:NLSDP-all}. A vector $\m q\in\mathbb{R}^d$ is called a feasible direction for problem~\eqref{eq:NLSDP-all} at $\m p^0$ if $\m p^0 + \varepsilon \m q$ is feasible for~\eqref{eq:NLSDP-all} for all sufficiently small $\varepsilon>0$.
\end{mydef}
The KKT conditions are of first order, i.e., they involve the gradient of the Lagrangian. The following definition states a second-order condition.
\begin{mydef}
\label{def:2ndordersuff}
Let $(\m p^*, \m \Lambda^*)$ be a KKT point of~\eqref{eq:NLSDP-all}. The second-order sufficiency condition holds for $\m p^*$ if for all feasible directions $\m q$ at $\m p^*$ satisfying $\nabla_{\m p}f(\m p^*)^{\top}\m q = 0$, it holds that $\m q^{\top}\nabla_{\m p}^{2}\mathcal{L}(\m p^*,\m \Lambda^*)\m q \geq \mu\|\m q\|^2$, for some $\mu>0$. 
\end{mydef}
The second-order sufficiency condition will be useful for the convergence of one of the algorithms to solve BMIs in the sequel. 
Sections~\ref{sec:relaxMIBMIs} and~\ref{sec:succapproxMIBMIs} develop algorithms for solving problems of the form~\eqref{eq:NLSDP-all} that include BMIs. These algorithms typically return vectors $\m \pi$ with non-integer, real entries. Based on the solutions produced by these algorithms,  Appendix~\ref{sec:slicing} details the procedure of actuator selection.  
\section{SDP Relaxations (SDP-R): A Lower Bound on~\eqref{eq:NLSDP-all}}\label{sec:relaxMIBMIs}
This section develops a solver for BMI problems based on SDP relaxation of the BMI constraint. To this end, we introduce an additional optimization variable $\m G  = \m \Pi  \m Z $. With this change of variables, $ \m \Pi  \m Z $ is replaced by $\m G $ and $\m G^{\top}$ replaces  $\mZ^{\top} \m\Pi$ in \eqref{eq:StabUI2-MIBMI}, while the constraint $\m G  = \m \Pi  \m Z $ is added to the problem. Effectively, we have pushed the bilinearity into a new constraint $\m G  = \m \Pi  \m Z $, which can actually be manipulated to much simpler constraints due to the diagonal structure of~$\m \Pi $. 

Specifically, $\m Z $ and $\m G $ are stacks of $N$ matrices
\begin{equation}
\mZ =\begin{bmatrix}
\mZ_1  \\ \vdots \\ \mZ_N  
\end{bmatrix}, \quad
\mG  =\begin{bmatrix}
\mG_1   \\ \vdots \\ \mG_N  
\end{bmatrix}
\label{eq:ZandGblock}
\end{equation} 
where $\mZ_i$ and $\mG_i$ ($i=1,\ldots,N$) are both in $\mathbb{R}^{n_{u_i}\times n_x}$. Due to the diagonal structure of $\m \Pi $, the constraint $\mG  = \m\Pi  \mZ $ is equivalent to
\begin{equation}
\mG_i  = \pi_i  \mZ_i , \quad i=1,\ldots,N.
\label{eq:bilin-equiv1}
\end{equation}
Denote the $(l,m)$ entries of matrices $\m Z_i $ and $\m G_i $ by $Z_{i,(l,m)} $ and $G_{i,(l,m)} $, respectively, where $l=1,\ldots,n_{u_i}$ and $m=1,\ldots,n_x$. Then,~\eqref{eq:bilin-equiv1} is equivalent to the constraint
\begin{align}
 G_{i,(l,m)}  = \pi_i    Z_{i,(l,m)}  , \quad & i=1,\ldots,N, \;  l=1,\ldots,n_{u_i},  \notag \\
 &  m=1,\ldots,n_x.
\label{eq:bilin-equiv2}
\end{align}
It follows that problem~\eqref{eq:NLSDP-all} is equivalent to
\begin{subequations}
				\label{eq:StabUI3-all}
				\begin{eqnarray}
				L= \minimize_{\mS, \mZ,  \zeta, \m\pi, \mG}& &\hspace{-0.41cm} \zeta   + \m \alpha_{\pi}^{\top} \m \pi \\
\subjectto 				& & \hspace{-0.651cm}\bmat{\mA \mS +\mS \mA^{\top}+\alpha\mS  \\-\mB_u  \mG  -\mG^{\top} \mB_u^{\top} &\mB_w  \\\mB_w^{\top} & -\alpha\eta \mI} \preceq \m {O} \label{eq:StabUI3-LMI}\\
%	& &   \hspace{-0.41cm}G_{i,(l,m)}   = \pi_i    Z_{i,(l,m)}  ,  i=1,\ldots,N, \; \\
%					&& \hspace{-0.41cm} l=1,\ldots,n_{u_i}, m=1,\ldots,n_x.  \label{eq:StabUI3-bilin}\\
					& &\hspace{-0.41cm} \eqref{eq:StabUI2-zetaS},\eqref{eq:StabUI2-coupl},\eqref{eq:StabUI2-box} \label{eq:StabUI3-rep}, \eqref{eq:bilin-equiv2}.
				\end{eqnarray}
			\end{subequations}		
The next step is to relax~\eqref{eq:bilin-equiv2} into an SDP constraint. To this end, define 
\begin{equation}
\mE = \begin{bmatrix}
0 & 0 & 0 \\
0 & 0 & 1 \\
0 & 1 & 0 \\
\end{bmatrix}, \quad \m e= \begin{bmatrix}
2\\ 0 \\ 0
\end{bmatrix}.
\end{equation}
The SDP relaxation of~\eqref{eq:StabUI3-all} is provided in the next proposition.
\begin{prop}
\label{prop:SDPrelax}
The following SDP is a relaxation of~\eqref{eq:StabUI3-all} and yields a lower bound on the optimal value of~\eqref{eq:NLSDP-all}
\begin{subequations}
				\label{eq:StabUI-sdprelax-all}
				\begin{eqnarray}
				\tilde{L}= \minimize_{\mS, \mZ,  \zeta, \m\pi, \mG, \mV}\hspace{1.2cm} & &\hspace{-1.8cm} (\eta+1)\zeta + \m \alpha_{\pi}^{\top} \m \pi \\
\subjectto 				\hspace{1.5cm}& & \nonumber\\
				&&\hspace{-4cm} \bmat{\mA \mS +\mS \mA^{\top}+\alpha\mS  \\-\mB_u  \mG  -\mG^{\top} \mB_u^{\top} &\mB_w  \\\mB_w^{\top} & -\alpha\eta \mI} \preceq \m {O} \label{eq:StabUI-sdprelax-LMI}\\
					& &  \hspace{-4cm} \mathbf{trace}\left( \mE \mV_{i,(l,m)}  \right) - \m e^{\top} \begin{bmatrix}
G_{i,(l,m)}   \\   Z_{i,(l,m)}   \\ \pi_i  
\end{bmatrix}=0   \label{eq:StabUI-sdprelax-trace} \\
					& &   \hspace{-4cm} \begin{bmatrix} 
\mV_{i,(l,m)}  &\left| \begin{matrix}
G_{i,(l,m)}   \\   Z_{i,(l,m)}   \\ \pi_i  
\end{matrix}\right.\\
\hline
\begin{matrix}
G_{i,(l,m)}   &  Z_{i,(l,m)}   & \pi_i  
\end{matrix} & 1
\end{bmatrix} \succeq \m O \label{eq:StabUI-sdprelax-psd} \\
& &  \hspace{-4cm}  \forall \; i=1,\ldots,N, \,  l=1,\ldots,n_{u_i}, \,  m=1,\ldots,n_x \notag \\
& &\hspace{-4cm}\eqref{eq:StabUI2-zetaS},\eqref{eq:StabUI2-coupl},\eqref{eq:StabUI2-box}
				\end{eqnarray}
			\end{subequations}		
where $\mV_{i,(l,m)} \in\mathbb{R}^{3 \times 3}$ are auxiliary optimization variables collected in $\mV$	 for all $i$, $l$, and $m$. 	
The optimal value of~\eqref{eq:StabUI-sdprelax-all} has the property that $\tilde{L}\leq L$. If in addition $\mathbf{rank}\left[\m V_{i,(l,m)}\right] = 1$ holds for all  $i$, $l$, and $m$ for the solution of~\eqref{eq:StabUI-sdprelax-all}, then  $\tilde{L} = L$. 
\end{prop}

\begin{IEEEproof}[Proof of Proposition~\ref{prop:SDPrelax}]
Introduce an auxiliary optimization variable
\begin{equation}
\m v_{i,(l,m)}  =\begin{bmatrix}
G_{i,(l,m)}   \\   Z_{i,(l,m)}   \\ \pi_i  
\end{bmatrix} \in \mathbb{R}^3
\label{eq:def31}
\end{equation}
With the previous definitions, it can easily be verified that
\begin{equation}
\pi_i    Z_{i,(l,m)}   - G_{i,(l,m)}   = \m v_{i,(l,m)}^{\top} \mE \m v_{i,(l,m)}  -\m e^{\top} \m v_{i,(l,m)} .
\end{equation}
A relaxation trick can be used at this point. In particular, introduce an additional optimization variable $\mV_{i,(l,m)} \in\mathbb{R}^{3 \times 3}$ and the constraint $\m V_{i,(l,m)}  = \m v_{i,(l,m)}  \m v_{i,(l,m)}^{\top}$. We have that 
\begin{align}
\m v_{i,(l,m)}^{\top} \mE \m v_{i,(l,m)}   & = \mathbf{trace}\left(\m v_{i,(l,m)}^{\top} \mE \m v_{i,(l,m)}  \right) \notag \\
& = \mathbf{trace}\left( \mE \m v_{i,(l,m)}   \m v_{i,(l,m)}^{\top}  \right) \notag \\
& = \mathbf{trace}\left( \mE \mV_{i,(l,m)}  \right). \label{eq:sdprelax-trace}
\end{align}
The previous development reveals that constraint~\eqref{eq:bilin-equiv2} is equivalent to the constraint $\mathbf{trace}\left( \mE \mV_{i,(l,m)}  \right) - \m e^{\top} \m v_{i,(l,m)}  =0$, which is linear in $\mV_{i,(l,m)} $ and  $\m v_{i,(l,m)}$, as long as the constraint $\m V_{i,(l,m)}  = \m v_{i,(l,m)}  \m v_{i,(l,m)}^{\top}$ is imposed, which is nonconvex.
The constraint $\m V_{i,(l,m)}  = \m v_{i,(l,m)}  \m v_{i,(l,m)}^{\top}$ is equivalent to
\begin{equation}
\begin{bmatrix}
\mV_{i,(l,m)}  & \m v_{i,(l,m)}  \\
\m v_{i,(l,m)}^{\top} & 1
\end{bmatrix} \succeq \m O, \; \mathbf{rank}(\m V_{i,(l,m)}  ) = 1.
\end{equation}

The rank constraint above is nonconvex, and by dropping it, we obtain the convex relaxation~\eqref{eq:StabUI-sdprelax-all} of~\eqref{eq:StabUI3-all}.
As a relaxation of~\eqref{eq:StabUI3-all}, its optimal value has the property that $\tilde{L}\leq L$. 
\end{IEEEproof} 
\normalcolor

Proposition~\ref{prop:SDPrelax} asserts that $\tilde{L}=L$ if $\mathbf{rank}\left[\m V_{i,(l,m)}\right] = 1$. Since the rank constraint is nonconvex, it is reasonable to consider surrogates for the rank in an effort to make the relaxation~\eqref{eq:StabUI-sdprelax-all}  tighter; one such surrogate is the nuclear norm of a matrix~\cite{recht2010guaranteed}. Thus, the constraint  $\normof{\mV_{i,(l,m)} }_{*}\leq 1$ can be added to promote smaller rank for ${\mV_{i,(l,m)} }$; the optimal value  of~\eqref{eq:StabUI-sdprelax-all} is impacted as follows. 
\begin{corol}
Let $\breve{L}$ be the optimal value of~\eqref{eq:StabUI-sdprelax-all} with the added constraint $\normof{\mV_{i,(l,m)} }_{*}\leq 1$. It holds that $\breve{L}\geq \tilde{L}$.
\label{corol:nucnorm}
\end{corol}
 
\begin{IEEEproof}[Proof of Corollary~\ref{corol:nucnorm}]
Adding the constraint restricts the feasible set of~\eqref{eq:StabUI-sdprelax-all}, yielding the stated relationship between the optimal values. 
\end{IEEEproof}

\normalcolor

\section{Convex Approximations: An Upper Bound on~\eqref{eq:NLSDP-all}}\label{sec:succapproxMIBMIs}

The common thread between the previous and the present section is to replace the nonconvex feasible set given by constraints~\eqref{eq:StabUI2-MIBMI}, \eqref{eq:StabUI2-zetaS},~\eqref{eq:StabUI2-coupl}, and~\eqref{eq:StabUI2-box} with convex sets. While the previous section relies on convex relaxations of the nonconvex feasible set, this section develops convex restrictions, i.e., replaces the nonconvex feasible set with a convex subset. The premise is to solve a series of optimization problems, in which the convex subset is improved. Thus, the algorithms in this section fall under the class of \emph{successive convex approximations} (SCAs). Two SCA algorithms are developed in this section.
Due to the convex restriction, the algorithms solve optimization problems that yield upper bounds for the optimal value $L$ of problem~\eqref{eq:NLSDP-all}. 

Because the SCA algorithms rely on forming convex subsets of the feasible nonconvex set, they must be initialized at interior points of the nonconvex feasible set. The next proposition asserts that such points indeed exist under Assumption~\ref{asmp:Feas}. 
\begin{prop}
\label{prop:SlaterRelax}
Under Assumption~\ref{asmp:Feas}, problem~\eqref{eq:NLSDP-all} is strictly feasible, i.e., it satisfies Slater's constraint qualification.  
\end{prop}
 
\begin{IEEEproof}[Proof of Proposition~\ref{prop:SlaterRelax}]
Consider a point $\m p^0$ that satisfies Assumption~\ref{asmp:Feas} (in particular, $\bm{\pi}^{0}=\bm{1}$ holds). Constraints~\eqref{eq:StabUI2-MIBMI}, \eqref{eq:StabUI2-zetaS},~\eqref{eq:StabUI2-coupl}, and~\eqref{eq:StabUI2-box} can be written in the form of a block diagonal matrix inequality~\eqref{eq:NLSDP-MI}. The implication is that $\mathcal{G}(\m p^0)$ is negative definite, i.e., all its eigenvalues are negative. By continuity of the eigenvalues as functions of the matrix elements~\cite[Appendix~D]{horn2012matrix}, there is a ball of sufficiently small radius around $\m p_0$ such that for all $\m p$ is this ball, the eigenvalues of $\mathcal{G}(\m p)$ remain negative. Any point within the ball satisfying $\bm{\pi}< \bm{1}$ together with the associated $\m S, \zeta, \m Z$ yields a stritly feasible point for constraints~\eqref{eq:StabUI2-MIBMI}, \eqref{eq:StabUI2-zetaS},~\eqref{eq:StabUI2-coupl}, and~\eqref{eq:StabUI2-box}.
\end{IEEEproof}
\normalcolor

\subsection{SCA using difference of convex functions (SCA-1)}
\label{subsec:SuccConvApprox}

The main idea is to replace~\eqref{eq:StabUI2-MIBMI} with a surrogate convex inequality constraint.  To this end, the left-hand side of~\eqref{eq:StabUI2-MIBMI} is replaced by a convex function in the variables $\mZ$, $\m \Pi$, which is denoted by $\mathcal{C} (\m \Pi , \mZ  ;\m\Pi _0, \mZ_0 )$, where $\m\Pi _0, \mZ_0 $ are given matrices to be specified later. This approach has been investigated in the context of BMIs for control problems with bilinearities arising in output feedback control problems; see~\cite{dinh2012combining}.
We first define the following linear function of $\m \Pi , \mZ$ with parameters  $\m\Pi _0, \mZ_0$
\begin{align}
& \mathcal{H}_{\mathrm{lin}} (\m\Pi ,\mZ ;\m\Pi _0, \mZ_0 ) = \notag \\
 & +\: \mB_u  \m\Pi_0  \m\Pi_0^{\top}  \mB_u^{\top}   - \mB_u  \m\Pi  \m\Pi_0^{\top}  \mB_u^{\top}  
 - \mB_u  \m\Pi_0  \m\Pi^{\top}  \mB_u^{\top}  \notag\\
 & +\: \mB_u  \m\Pi_0  \mZ_0^{j} - \mB_u  \m\Pi  \mZ_0^{j}  -  \mB_u  \m\Pi_0  \mZ^{j} \notag\\
 & +\: \mZ_0^{\top} \m\Pi_0  \mB_u^{\top} -  \mZ_0^{\top} \m\Pi  \mB_u^{\top}  - \mZ^{\top} \m\Pi_0  \mB_u^{\top}  \notag \\
 & +\: \mZ_0^{\top} \mZ_0 - \mZ_0^{\top} \mZ - \mZ^{\top} \mZ_0.
  \label{eq:HLinDef}
\end{align}
The function $\mathcal{C} (\m \Pi , \mZ  ;\m\Pi _0, \mZ_0 )$ is given by
\begin{multline}
\mathcal{C} (\cdot )\hspace{-0.02cm}=\hspace{-0.1cm}\bmat{\mA \mS +\mS \mA^{\top}+\alpha\mS  \\ +\frac{1}{2}  \left( \mB_u  \m\Pi  - \mZ^{\top}\right) \left(\mB_u  \m\Pi  - \mZ^{\top}\right)^{\top} \\ +\frac{1}{2}  \mathcal{H}_{\mathrm{lin}} (\m\Pi ,\mZ ;\m\Pi _0, \mZ_0 ) &\mB_w  \\\mB_w^{\top} & -\alpha\eta \mI}.
\label{eq:CDef}
\end{multline}
The following proposition asserts that $\mathcal{C} (\m \Pi , \mZ  ;\m\Pi _0, \mZ_0 )$ is a convex function that upper bounds the left-hand side of~\eqref{eq:StabUI2-MIBMI}.
\begin{prop}
\label{prop:CUppBd}
It holds for all $\m \Pi , \mZ$ and $\m\Pi _0, \mZ_0$ that 
\begin{equation}
 \hspace{-0.11cm}\bmat{\mA \mS +\mS \mA^{\top}+\alpha\mS  \\-\mB_u  \m\Pi  \mZ  -\mZ^{\top} \m\Pi  \mB_u^{\top} &\mB_w  \\\mB_w^{\top} & -\alpha\eta \mI} \preceq \mathcal{C} (\m \Pi , \mZ  ;\m\Pi _0, \mZ_0 ),
 \label{eq:BMIuppbd}
\end{equation}
where $$\small \mathcal{C} (\m \Pi , \mZ  ;\m\Pi _0, \mZ_0 ) = \bmat{\mA \mS +\mS \mA^{\top}+\alpha\mS  \\ +\frac{1}{2}  \left( \mB_u  \m\Pi  - \mZ^{\top}\right) \left(\mB_u  \m\Pi  - \mZ^{\top}\right)^{\top} \\ +\frac{1}{2}  \mathcal{H}_{\mathrm{lin}} (\m\Pi ,\mZ ;\m\Pi _0, \mZ_0 ) &\mB_w  \\\mB_w^{\top} & -\alpha\eta \mI}$$ is convex in $\m \Pi , \mZ$. 
\end{prop}
{The proof of Proposition~\ref{prop:CUppBd} is included in Appendix~\ref{app:allproofs}}. Given this result, convex approximation of the BMI is obtained by replacing constraint~\eqref{eq:StabUI2-MIBMI} with the convex constraint
$\mathcal{C} (\m \Pi , \mZ  ;\m\Pi _0, \mZ_0 )\preceq \mO$.
The resulting problem has a restricted feasible set due to~\eqref{eq:BMIuppbd}. 
Although $ \mathcal{C} (\m \Pi , \mZ  ;\m\Pi _0, \mZ_0 )$  is a convex function in $\m\Pi $ and $\mZ $, it is not linear in $\m\Pi $ and $\mZ $.
Therefore, when we replace~\eqref{eq:StabUI2-MIBMI}  by the constraint  $\mathcal{C} (\m \Pi , \mZ  ;\m\Pi _0, \mZ_0 )\preceq \mO$, a convex constraint is obtained, but not an LMI. Fortunately, the constraint  $\mathcal{C} (\m \Pi , \mZ  ;\m\Pi _0, \mZ_0 )\preceq \mO$ can be equivalently written as an LMI as follows.
\begin{lem}
\label{lem:CsLMI}
It holds that
\begin{multline}
 \mathcal{C} (\m \Pi , \mZ  ;\m\Pi _0, \mZ_0 ) \preceq \mO \Longleftrightarrow  \mathcal{C}_{s} (\m \Pi , \mZ  ;\m\Pi _0, \mZ_0 ) =\\
\small   \bmat{\mA \mS +\mS \mA^{\top}+\alpha\mS   \\ +\frac{1}{2}  \mathcal{H}_{\mathrm{lin}} (\m\Pi ,\mZ ;\m\Pi _0, \mZ_0 ) 
   & \frac{1}{\sqrt{2}} \left( \mB_u  \m\Pi  - \mZ^{\top}\right) & \mB_w  \\
   \frac{1}{\sqrt{2}}\left(\mB_u  \m\Pi  - \mZ^{\top}\right)^{\top} & -\mI & \mO  \\
   \mB_w^{\top} & \mO  & -\alpha\eta \mI}  \preceq \mO.
   \label{eq:BMIUppBd-LMI}
\end{multline}
\end{lem}
 
\begin{IEEEproof}[Proof of Lemma~\ref{lem:CsLMI}]
Applying the Schur complement to $\mathcal{C} (\m \Pi , \mZ  ;\m\Pi _0, \mZ_0 )\preceq \mO$ yields the LMI $ \mathcal{C}_{s} (\m \Pi , \mZ  ;\m\Pi _0, \mZ_0 )$.
\end{IEEEproof}
\normalcolor

To summarize, the convex approximation to~\eqref{eq:StabUI2-all} at $\m\Pi _0, \mZ_0 $ is formed by replacing  the integer constraints by the box constraints~\eqref{eq:StabUI2-box}, and the BMI~\eqref{eq:StabUI2-MIBMI} by the constraint the LMI in~\eqref{eq:BMIUppBd-LMI}.
This problem is stated as follows:
\begin{subequations}
				\label{eq:StabUI-DC-all}
				\begin{eqnarray}
				\hat{L}= \minimize_{\mS, \mZ,  \zeta, \m\pi} && (\eta+1)\zeta  + \m \alpha_{\pi}^{\top} \m \pi \label{eq:StabUI-DC-obj}\\
\subjectto 			&&	\eqref{eq:StabUI2-zetaS},\eqref{eq:StabUI2-coupl},\eqref{eq:StabUI2-box},\eqref{eq:BMIUppBd-LMI}.		
				\end{eqnarray}
			\end{subequations}		
Problem~\eqref{eq:StabUI-DC-all} is an SDP with optimal value denoted by $\hat{L}$, whose relationship with $L$ is as follows. 
\begin{corol}
\label{corol:SCA1optval}
The optimal value of the convex approximation~\eqref{eq:StabUI-DC-all} for all $\m\Pi _0, \mZ_0$ is an upper bound on the optimal value of~\eqref{eq:NLSDP-all}, that is, $L\leq \hat{L}$.
\end{corol}
 
\begin{IEEEproof}[Proof of Corollary~\ref{corol:SCA1optval}]
Due to~\eqref{eq:BMIuppbd} and~\eqref{eq:BMIUppBd-LMI}, problem~\eqref{eq:StabUI-DC-all} has a restricted feasible set with respect to problem~\eqref{eq:NLSDP-all}. 
\end{IEEEproof}
\normalcolor
The convex approximation~\eqref{eq:StabUI-DC-all} depends on the point $\m\Pi _0, \mZ_0 $, and can be successively improved. The main idea is to solve a sequence of convex approximations given by~\eqref{eq:StabUI-DC-all}, where the values of $\m\Pi _0, \mZ_0 $ for  the next approximating problem are given by the solution of the previous problem. 

Let $k=1,2,\ldots$ denote the index of the convex approximation to be solved, and let $\mS_k, \zeta_k, \m\Pi_k, \mZ_k$ denote its solution. The $k$-th problem is obtained by {  adding a strictly convex regularizer to the objective~\eqref{eq:StabUI-DC-obj}, which ensures that the problem has a unique solution.}  The $k$-th problem is thus
\begin{subequations}
				\label{eq:StabUI-DC2-all}
				\begin{align}
				\hat{L}^{(1)}_{k}= \minimize_{\{\mS, \mZ,\zeta,\m\pi\}} &\;\;\;(\eta+1)\zeta + \m \alpha_{\pi}^{\top}\m \pi  +  \rho J_{k}\label{eq:StabUI-DC2-obj}\\
%\subjectto 	&&\nonumber \\
\subjectto&\;\;\; \mathcal{C}_s (\m \Pi , \mZ  ;\m\Pi_{k-1}, \mZ_{k-1} ) \preceq \m O \\
%& \;\;\;\bmat{-\mS  & \mO & \mS \mC_z^{ \top}\\
%	\mO & -\mI & \mD_{wz}^{\top}\\
%	\mC_z \mS & \mD_{wz}  &-\zeta \mI} \preceq \m {O} \\
	&\;\;\; \eqref{eq:StabUI2-zetaS},\m H \m\pi \leq \m h, \; \;0 \leq \m \pi \leq 1,	
%& \eqref{eq:StabUI2-zetaS},
				\end{align}
			\end{subequations}		
where $J_{k}= \|\zeta - \zeta_{k-1}\|_2^2 +  \|\mS - \mS_{k-1}\|_F^2 + \|\mZ - \mZ_{k-1}\|_F^2 + \|\m\Pi - \m\Pi_{k-1}\|_F^2$; the linearization point is given by $\m\Pi _0=\m\Pi_{k-1}$, $\mZ_0 =\mZ_{k-1}$; $\rho$ is the weight of the quadratic regularizers. For $k=1$, the point $\m\Pi _0, \m Z_0$ can be selected as any interior point of~\eqref{eq:NLSDP-all}; such is guaranteed to exist due to Proposition~\ref{prop:SlaterRelax}. {  Note that the regularization term $\rho J_k$ penalizes the difference between the new solution and the previous. Upon algorithm convergence, the two successive solutions should be close to each other, which means that at optimality, the entire term $\rho J_k$ should be close to zero.}

Notice that for every $k$, problem~\eqref{eq:StabUI-DC2-all} has the form of~\eqref{eq:NLSDP-all}, but the objective is a strictly convex quadratic, and the constraint function is convex. The convergence is established in the following proposition.
\begin{prop}
\label{prop:DCconv}
Let $\m p_k, \m \Lambda_k$ denote a KKT point of~\eqref{eq:StabUI-DC2-all}. Suppose that the feasible set of~\eqref{eq:NLSDP-all} is bounded, and that  the following hold for problem~\eqref{eq:StabUI-DC2-all} for $k=1,2,3,...$
\begin{itemize}\setlength{\itemindent}{.1in}
\item[i)] Slater's constraint qualification holds.
\item[ii)] The Lagrange multiplier  $\m \Lambda_k$ is locally unique.
\item[iii)] Strict complementarity holds for the KKT point.
\item[iv)] The second-order sufficiency condition holds for the KKT point. 
\end{itemize}
Then, the following are concluded:
\begin{itemize}\setlength{\itemindent}{.1in}
\item[a)] It holds that $f(\m p_k)\geq L$ and $L_{k}^{(1)}\geq L$ for $k=1,2,3,...$
\item[b)] The sequence $\{f(\m p_k)\}_{k=1}^\infty$ is monotone decreasing, and converges to a limit $\hat{f}^{(1)}\geq L$.
\item[c)] Every limit point of the sequence $\{\m p_k, \m \Lambda_k\}_{k=1}^\infty$ is a KKT point of~\eqref{eq:NLSDP-all}. If the set of KKT points of~\eqref{eq:NLSDP-all} is finite, then the entire sequence  $\{\m p_k, \m \Lambda_k\}_{k=1}^\infty$ converges to a KKT point of~\eqref{eq:NLSDP-all}.
\end{itemize}
\end{prop}
{The proof of Proposition~\ref{prop:DCconv} is included in Appendix~\ref{app:allproofs}}.  Albeit some of the conditions of the previous proposition may be hard to verify in practice, we encountered no case where the SCA algorithm did not converge. In particular, we tested the algorithm on a variety of dynamic systems with varying sizes and conditions in Section~\ref{sec:numtests}.

\subsection{Parametric SCA (SCA-2)}\label{subsec:paramSCA}
In this section, we depart from the  difference of two convex functions approach used in the previous SCA, and use another approach to obtain an upper bound on the bilinear terms. The developments in this section follow the spirit of the methods presented in~\cite{LeeHu2016}, where the authors investigate a new approach to solve BMIs that are often encountered in output feedback control problems. 

First, let $\mathcal{F}_1(\m p)$ denote the left-hand side of~\eqref{eq:StabUI2-MIBMI}.
Given $\m\Pi_{k}$ and $\mZ_{k}$, define $\Delta \m\Pi = \m\Pi - \m\Pi_{k} $ and $\Delta \mZ = \mZ - \mZ_{k}$. For any $\boldsymbol{Q} \in \mathbb{S}^{n_x}_{++}$, define further the following function:
\begin{align}
\mathcal{\m K}_1 (\m p;\m p_{k},\m Q) = 
	&\bmat{-\mB_{u} \m\Pi_{k} \mZ_{k} - \mZ_{k} ^{\top}\m\Pi_{k}\mB_{u}^{\top} & \mB_{w} \\  \mB_{w} ^{\top} & -\alpha\eta \mI}\notag\\ 	
	&\mspace{-100mu}+\bmat{\mA\mS+\mS\mA^{\top}+\alpha\mS -\mB_{u}\m\Pi_{k}\Delta \mZ & \\ - \Delta \mZ ^{\top}\m\Pi _{k}\mB_{u} ^{\top} - \mB_{u}\Delta\m\Pi \mZ_{k} - \mZ_{k} ^{\top}\Delta \m\Pi \mB_{u} ^{\top} & \mO \\ \mO & \mO} \notag \\ 
	&\mspace{-100mu}+\bmat{\mB_{u} \Delta \m\Pi \mQ \Delta \m\Pi \mB_{u}^{\top}+ \Delta \mZ ^{\top} \mQ^{-1}\Delta \mZ & \mO \\ \mO & \mO}.
	\label{eq:Kscadef}
\end{align}

Similar to Proposition~\ref{prop:CUppBd}, an upper bound on  $\mathcal{\m F}_1(\m p)$ is provided by the next proposition. 
\begin{prop}
\label{prop:K1UppBd}
It holds for all $\m p$, $\m\Pi_{k}, \mZ_{k}$ and $\m Q \in \mathbb{S}^{n_u}_{++}$ that
\begin{equation}
\mathcal{\m F}_1(\m p)  \leq \mathcal{\m K}_1 (\m p;\m p_{k},\m Q).
\label{eq:K1UppBd}
\end{equation}
\end{prop}
{The proof of Proposition~\ref{prop:K1UppBd} is included in Appendix~\ref{app:allproofs}}. The previous proposition suggests that constraint~\eqref{eq:StabUI2-MIBMI} can be replaced by $\mathcal{\m K}_1 (\m p;\m p_{k},\m Q) \preceq \m O$. There are two challenges to be addressed though. First, although $\mathcal{\m K}_1 (\m p;\m p_{k},\m Q)$ is a convex function of $\m p$, it is not linear, and thus constraint $\mathcal{\m K}_1 (\m p;\m p_{k},\m Q) \preceq \m O$ is not an LMI. Second, although $\m Q$ can remain constant, the approximation can be tightened if $\m Q$ is allowed to be an optimization variable. The former challenge is addressed by Lemma~\ref{lem:KLMI}, which is analogous to Lemma~\ref{lem:CsLMI}.
\begin{lem}
\label{lem:KLMI}
Constraint $\mathcal{\m K}_1 (\m p;\m p_{k},\m Q) \preceq \m O$ is equivalent to
\begin{equation}
\small \mathcal{\m K}(\m p;\m p_{k},\m Q) = \bmat {\m \Omega(\m p; \m p_k) & \mB_{w} & \mB_{u}\Delta \m\Pi & \Delta \mZ^{\top}\\
		\mB_{w}^{\top} & -\alpha\eta \mI & \mO  & \mO \\
		\Delta \m\Pi \mB_{u}^{\top} & \mO  & -\mQ^{-1} & \mO \\
		\Delta \mZ & \mO  & \mO  & -\mQ} \preceq \m O,
\end{equation}
\begin{align*}
&\m \Omega(\m p; \m p_k) = -\mB_{u} \m\Pi_{k} \mZ_{k} - \mZ_{k} ^{\top}\m\Pi_{k}\mB_{u}^{\top} 
 +\mA\mS+\mS\mA^{\top}+\alpha\mS \\
& \,\,-\mB_{u}\m\Pi_{k}\Delta \mZ
- \Delta \mZ ^{\top}\m\Pi _{k}\mB_{u} ^{\top}- \mB_{u}\Delta\m\Pi \mZ_{k} - \mZ_{k} ^{\top}\Delta \m\Pi \mB_{u} ^{\top}.
\end{align*}
\end{lem}
{\begin{IEEEproof}[Proof of Lemma~\ref{lem:KLMI}]
Use the Schur complement. 
\end{IEEEproof}}
When $\m Q$ is an optimization variable, function $\mathcal{\m K}(\m p;\m p_{k},\m Q)$ is not convex in $\m p$ and $\m Q$. An upper bound of $\mathcal{\m K}(\m p;\m p_{k},\m Q)$ that is linear in $\m p$ and $\m Q$ is given in Lemma~\ref{lem:KsUppBd}. The following lemma gives a particular matrix property that becomes the foundation for Lemma~\ref{lem:KsUppBd}. 
\begin{lem}
	\label{lem:MinQinvUppBd}
	Let $\mathcal{\mQ}(\m x):\mathbb{R}^{n}\rightarrow\mathbb{S}^{m}_{++}$ be a mapping defined as $\mathcal{\mQ}(\m x) = \sum_{i=1}^{n} x_i\mathcal{\mQ}_i$ where $\mathcal{\mQ}_i\in\mathbb{S}^{m}$. The following inequality holds, where the right-hand side is the linearization of $-\mathcal{\mQ}(\m x)^{-1}$ around $\m x_k$:
	\begin{equation}
	-\mathcal{\mQ}(\m x)^{-1} \preceq -2\mathcal{\mQ}(\m x_k)^{-1}+\mathcal{\mQ}(\m x_k)^{-1}\mathcal{\mQ}(\m x)\mathcal{\mQ}(\m x_k)^{-1}.
	\label{eq:MinQinvUppBd}
	\end{equation}
\end{lem}

\begin{lem}
	\label{lem:KsUppBd}
	It holds for all $\m p$, $\m Q\in \mathbb{S}^{n_u}_{++}$,  $\m\Pi_{k}, \mZ_{k},$ and $\m Q_{k} \in \mathbb{S}^{n_u}_{++}$ that
	\begin{equation}
	\mathcal{\m K}(\m p;\m p_{k},\m Q) \preceq \mathcal{\m K}_s (\m p, \m Q;\m p_{k},\m Q_k)
	\label{eq:KsUppBd}
	\end{equation}
	where $\mathcal{\m K}_s (\m p,\m Q;\m p_{k},\m Q_k) =$
	\begin{equation}
	\bmat {\m \Omega(\m p; \m p_k) & \m B_w & \mB_{u}\Delta \m\Pi\mQ_{k} & \Delta \mZ^{\top}\\
		\mB_{w}^{\top} & -\alpha\eta \mI & \mO  & \mO \\
		\mQ_{k}\Delta \m\Pi \mB_{u}^{\top} & \mO  & -2\mQ_{k}+\mQ & \mO \\
		\Delta \mZ & \mO  & \mO  & -\mQ}.
	\label{eq:Ksdef}
	\end{equation} 
\end{lem}
{The proofs of Lemmas~\ref{lem:MinQinvUppBd} and~\ref{lem:KsUppBd} are included in Appendix~\ref{app:allproofs}}. Given these results, the constraint $\mathcal{\m K}_s (\m p, \m Q;\m p_{k},\m Q_k)\preceq \m O$ yields a restricted feasible set relative to constraint~\eqref{eq:StabUI-MIBMI}. Similarly to Section~\ref{subsec:SuccConvApprox}, $k=1,2,3,...$ is the index of the optimization problem to be solved, and  $\m p_k, \mQ_k$ denotes its solution. The $k$-th problem is an SDP and is stated as follows.
\begin{subequations}~\label{eq:SecondCA}
	\begin{align}
		\hat{L}^{(2)}_{k} =	\underset{\{\mS, \mZ,  \zeta, \m\pi, \mQ\}}{\minimize} & \;\;\;(\eta+1)\zeta  + \m \alpha_{\pi}^{\top}\m \pi  + \rho J_{k}\\
		\subjectto 
		&\;\;\;\mathcal{\m K}_s (\m p, \m Q;\m p_{k-1},\m Q_{k-1}) \preceq \m O\\
%		&\;\;\;	\bmat{-\mS  & \mO & \mS \mC_z^{ \top}\\
%			\mO & -\mI & \mD_{wz}^{\top}\\
%			\mC_z \mS & \mD_{wz}  &-\zeta \mI} \preceq \m {O} \\
&\hspace{-1cm} c_{1}\m I \preceq \mQ \preceq c_{2}\m I,\,\, -2\mQ_{k-1}+\m Q\preceq -c_{3}\m I~\label{eq:SecondCA-auxcon}\\
		&\;\;\; \eqref{eq:StabUI2-zetaS}, \; \m H \m\pi \leq \m h, \; \;\bm{0} \leq \m \pi \leq \bm{1},	
	\end{align}
\end{subequations}
where $\rho$, $c_1$, $c_2$, and $c_3$ are positive constants, and $J_{k}$ is the same regularizer as the one in~\eqref{eq:StabUI-DC2-all}. Constraint~\eqref{eq:SecondCA-auxcon} guarantees that $\m Q$ is positive definite, sequence $\{\m Q_k\}_{k=1}^{\infty}$ is bounded, and that $ -2\mQ_{k}+\m Q$, which appears as a diagonal block in~\eqref{eq:Ksdef} is negative definite for all $k$. Similar to the first convex approximation, the above problem can be initialized by letting $\{\mS_0, \mZ_0,  \zeta_0, \m\pi_0\}$ be any interior point of~\eqref{eq:NLSDP-all} and $\m Q_0=\m I$. The algorithm convergence is characterized by the following proposition.
\begin{prop}
\label{prop:paramSCAconv}
Assume that the MFCQ holds for every feasible point of~\eqref{eq:NLSDP-all} and that the sequence $\{\m p_k\}_{k=1}^\infty$ is bounded. 
Then, the following are concluded:
\begin{itemize}\setlength{\itemindent}{.1in}
\item[a)] It holds that $f(\m p_k)\geq L$ and $L_{k}^{(2)}\geq L$ for $k=1,2,\ldots$
\item[b)] The sequence $\{f(\m p_k)\}_{k=1}^\infty$ is monotone decreasing, and converges to a limit $\hat{f}^{(2)}\geq L$.
\item[c)] Every limit point of $\{\m p_k\}_{k=1}^\infty$ is a stationary point of~\eqref{eq:NLSDP-all}.
\end{itemize}
\end{prop}
{The proof of Proposition~\ref{prop:paramSCAconv} is included in Appendix~\ref{app:allproofs}}. Algorithm~\ref{algorithm:SCA} in Appendix~\ref{app:scalgo} provides the option to implement one of the two developed convex approximations [cf.~\eqref{eq:StabUI-DC2-all} and~\eqref{eq:SecondCA}] sequentially until a maximum number of iterations ($\mathrm{MaxIter}$) or a stopping criterion defined by a tolerance ($\mathrm{tol}$) are met. The next section compares the two approximations in terms of computational effort and their convergence claims. 

\subsection{Comparing the  SCAs and Recovering the Integer Solutions}
\label{subsec:compareSCA}

The first convex approximation is simpler to implement and involves a smaller number of SDP constraints and variables; see the difference in dimensions between constraints $\mathcal{\m K}_s (\m p, \m Q;\m z_{k},\m Q_{k}) \preceq \m O$ and $\mathcal{C}_s (\m \Pi , \mZ;\m\Pi _k, \mZ_k ) \preceq \m O$. In addition, constraint~\eqref{eq:SecondCA-auxcon} is added, and an extra variable $\m Q$ is needed in~\eqref{eq:SecondCA}. %Hence, it is expected that the first convex approximation~\eqref{eq:StabUI-DC2-all} requires less computational time. 
Both methods rely on constructing a series of feasible sets that are subsets of the original nonconvex feasible set in~\eqref{eq:NLSDP-all}.  Each produces a sequence of decreasing objective values $\{f(\m p_k)\}_{k=1}^{\infty}$, yielding upper bounds on the optimal value of~\eqref{eq:NLSDP-all}.

It is also worth noting that the first method requires a constraint qualification and additional assumptions on the KKT point to hold for each convex approximation problem $k$. Slater's constraint qualification is also an assumption in one of the earliest SCA methods for nonlinear programming~\cite{marks1978general}. On the other hand, the second method requires only the MFCQ to hold for the original nonconvex problem~\eqref{eq:NLSDP-all}. Both methods have a boundedness assumption; the first method requires the feasible set of~\eqref{eq:NLSDP-all} to be bounded, the second method only the resulting sequence to be bounded. The boundedness assumption respectively guarantees the existence of at least one limit point of $\{\m p_k\}_{k=1}^{\infty}$. Both methods enjoy the property that every limit point of $\{\m p_k\}_{k=1}^{\infty}$ is a stationary point of~\eqref{eq:NLSDP-all}. 

\begin{rem}[Existence of Local Minima]
The stationarity is a necessary condition for local optimality (cf. Lemma~\ref{lem:1stOrderOpt}). It is thus not guaranteed that the stationary point is indeed locally optimal. In view of the fact that the methods attempt to solve a nonconvex problem, such convergence result is to be expected. It is worth asking whether the resulting limit point is indeed locally optimal. Sufficient conditions for local optimality of stationary points of nonlinear SDPs have been derived in the literature; see for example~\cite[Theorem 9]{shapiro1997first}. {  Note also that the stationary points that SCA-1 and 2 converge to depend in general on the initial linearization point $\m p_0$.}
	\end{rem} 
%\begin{rem}[Recovering the Optimal Solution of~\eqref{eq:NLSDP-all}]~\label{rem:recoveringBMIsol}
%A comparison with the SDP relaxation of Section~\ref{sec:relaxMIBMIs} is in order. The method of Section~\ref{sec:relaxMIBMIs} produces a lower bound $\tilde{L}$ for the optimal value $L$ of~\eqref{eq:NLSDP-all}. The methods of this section derive upper bounds ($\hat{L}^{(1)}$ or $\hat{L}^{(2)}$) of $L$. Therefore, if the difference between the upper and lower bounds is small, we can assert that we have approximately solved nonconvex problem~\eqref{eq:NLSDP-all}. 
%	\end{rem} 
%\subsection{Recovering the Integer Solutions}
The solutions obtained from~\eqref{eq:StabUI-sdprelax-all}, \eqref{eq:StabUI-DC2-all}, and~\eqref{eq:SecondCA} produce a noninteger solution for the actuator selection problem. Since the objective is to determine a binary selection for the actuators, we present in this section a simple \textit{slicing routine} that returns a binary selection given the solutions to the optimization problems in Sections~\ref{sec:relaxMIBMIs} and \ref{sec:succapproxMIBMIs}. The algorithm is included and discussed in Appendix~\ref{sec:slicing}.

\section{SaA Selection via MISDP and the Big-M Method}~\label{sec:bigM}
This section develops an alternative method for solving the optimization problem~\eqref{eq:StabUI-all}. This alternative can also be applied to other time-varying SaA selection problems with the control and estimation  metrics and formulations in Tables~\ref{tab:Controllers} and~\ref{tab:ObserverDesign} in Appendix~\ref{app:cont}.

As discussed in the previous sections, the mixed-integer bilinear term $\m B_u \m \Pi \m Z+\m Z^{\top}\m \Pi \mB_u^{\top}$ renders the problem nonconvex. An alternative to solving the  convex relaxations or approximations is to simply apply the Big-M method on the bilinear term. This technique is quite general~\cite{grossmann2002review} and has been used before in the context of multi-vehicle path planning~\cite{5717988}, and more recently for actuator and sensor allocation in linear systems with Gaussian disturbances and Kalman filtering~\cite{Taylor2017}. 

In order to state the Big-M method, we will use the block matrices defined in~\eqref{eq:ZandGblock}--\eqref{eq:bilin-equiv2}. In particular, notice that due to the binary nature of $\pi_i$, constraint~\eqref{eq:bilin-equiv2} can be equivalently written for all $l,m$  as
\begin{equation}
G_{i,(l,m)}=\begin{cases} 
\pi_i Z_{i,(l,m)}, &\text{if~}\pi_i=1 \\
0, &\text{if~}\pi_i=0.
\end{cases} 
\label{eq:bilin-equiv01}
\end{equation}

By introducing a sufficiently large constant $M$, it is shown that the previous constraint can be equivalently written as 
\begin{equation}
|G_{i,(l,m)}-Z_{i,(l,m)}| \leq M (1-\pi_i), |G_{i,(l,m)}| \leq M\pi_i
\label{eq:bilin-equiv-bigM}
\end{equation}
\begin{lem}
	\label{lem:bigMconstr}
	Under the constraint $\pi_{i}\in\{0,1\}$ for all $i$ and for sufficiently large $M$, then any $\m Z$ and $\m G$ satisfying~\eqref{eq:bilin-equiv01} also satisfy~\eqref{eq:bilin-equiv-bigM}, and vice versa. 
\end{lem}
\begin{IEEEproof}[Proof of Lemma~\ref{lem:bigMconstr}]
	Suppose that $\pi_i=1$. Then, both~\eqref{eq:bilin-equiv01} and~\eqref{eq:bilin-equiv-bigM} are equivalent to $G_{i,(l,m)}=Z_{i,(l,m)}$ for all $(l,m)$. Suppose now that $\pi_i=0$. Then, both~\eqref{eq:bilin-equiv01} and~\eqref{eq:bilin-equiv-bigM} are equivalent to $G_{i,(l,m)}=0$ for all $(l,m)$ and allow $Z_{i,(l,m)}$ to be free. Notice that the proof indicates how big the constant $M$ should be. In particular, it must be larger than the maximum entry of matrix $Z$ that is the optimal solution of~\eqref{eq:StabUI-all}. In practice, a very large number is selected. 
\end{IEEEproof}

The advantage of~\eqref{eq:bilin-equiv-bigM} is that is converts the bilinear constraint~\eqref{eq:bilin-equiv2} to a constraint \emph{linear} in $\mZ$, $\mG$ and $\m \Pi$. Therefore, the actuator selection problem for $\mathcal{L}_{\infty}$ control is written as
\begin{subequations}
	\label{eq:StabUI2-BigMall}
	\begin{align}
		B_M^* = 				\minimize_{\{\mS,\mZ, \zeta, \m\pi, \mG\}} & \;\;\; (\eta+1)\zeta +   \alpha_{\pi}^{\top} \m \pi \\
		\subjectto 		& \bmat{\mA\mS+\mS\mA^{\top}+\alpha\mS \\-\mB_u  \mG -\mG^{\top} \mB_u^{\top}&\mB_w \\\mB_w^{\top} & -\alpha\eta \mI} \preceq \m {O} \\
		& \bmat{-\mS & \mO & \mS\mC_z^{\top}\\
			\mO & -\mI & \mD_{wz}^{\top}\\
			\mC_z\mS& \mD_{wz} &-\zeta\mI} \preceq \m {O}\\
		& \m H \m\pi \leq \m h  ,\;\;\m\pi  \in \{0,1\}^N \\
		& \m \Theta_1 \m G+ \m \Theta_2 \m Z \leq \m L_M(\m \pi).
	\end{align}
\end{subequations}		
The last constraint in~\eqref{eq:StabUI2-BigMall} represents the vectorization of the Big-M constraints in~\eqref{eq:bilin-equiv-bigM}, where $\m \Theta_1$, $\m \Theta_2$ and $\mL (\m \pi)$ are suitable matrices, and the latter is a linear mapping in $\m \pi$ that depends on $M$. The optimization still includes the integrality constraints on $\m \pi$, and hence it is a mixed-integer semidefinite program (MISDP). The key point is that it is equivalent to~\eqref{eq:StabUI-all}.
\begin{prop}
	\label{prop:bigMequiv}
	For sufficiently big $M$, problems~\eqref{eq:StabUI2-BigMall} and~\eqref{eq:StabUI-all} are equivalent, and thus, have equal optimal values, i.e., $f^*=B_M^*$. 
\end{prop}
\begin{IEEEproof}[Proof of Proposition~\ref{prop:bigMequiv}]
	Introduce the change of variables $\m G = \m \Pi \mZ$ in~\eqref{eq:StabUI-all}. The resulting problem and~\eqref{eq:StabUI2-BigMall} have the same feasible sets due to Lemma~\ref{lem:bigMconstr}.
\end{IEEEproof}

The widely used convex optimization modeling toolbox YALMIP incorporates modeling of mixed-integer convex programs \cite{yalmipMI} and interface corresponding general-purpose SDP solvers combined with implementations of branch-and-bound (BB) methods. The BB method is essentially a smart and efficient way to run an exhaustive search over all possible $2^{N}$ combinations of the binary variables. At most $2^{N}$ SDPs are then solved in the worst case run of a BB method. However, the empirical complexity of BB algorithms is much smaller than the worst-case one. Thus, such off-the-shelf solvers can be applied to~\eqref{eq:StabUI2-BigMall}. For this purpose, we compare the performance of YALMIP's BB algorithm with the developed relaxations and approximations (SDP-R, SCA-1, SCA-2) in Section~\ref{sec:numtests}. 

\begin{rem}[Computational Complexity]\label{rem:CTIME}
	Primal-dual interior-point methods for SDPs have a worst-case complexity estimate of $\mathcal{O}\left(m^{2.75}L^{1.5} \right)$, where $m$ is the number of variables (a function of $N,n_x,n_u,n_z$) and $L$ is the number of constraints\cite{boyd1994linear}. {In various problems arising in control systems studies, it is shown that the complexity estimate is closer to $\mathcal{O}\left(m^{2.1}L^{1.2} \right)$ which is significantly smaller than the worst-case estimate $\mathcal{O}\left(m^{2.75}L^{1.5} \right)$~\cite{boyd1994linear}.}  Solving the SCAs involves iteratively obtaining a converging solution to the BMIs, and hence it is difficult to obtain an upper bound on the number of iterations and thus perform any comparison of SCAs with the MISDP~\eqref{eq:StabUI2-BigMall}.  But the worst-case complexity of~\eqref{eq:StabUI2-BigMall} is $\mathcal{O}\left(2^N m^{2.75}L^{1.5} \right)$ (notice the exponential factor). As for the SDP relaxation, the computational complexity is only that of an SDP, hence it scales better than MISDPs or the SCAs.
\end{rem}

\begin{rem}
	Replacing the integrality constraint on $\m \pi$ with the box constraint in~\eqref{eq:StabUI2-BigMall} yields an SDP that can be solved using classical SDP solvers. To obtain the binary actuator selection, Algorithm~\ref{algorithm:Slicing2} can be implemented. This can significantly reduce the computational time.
\end{rem}

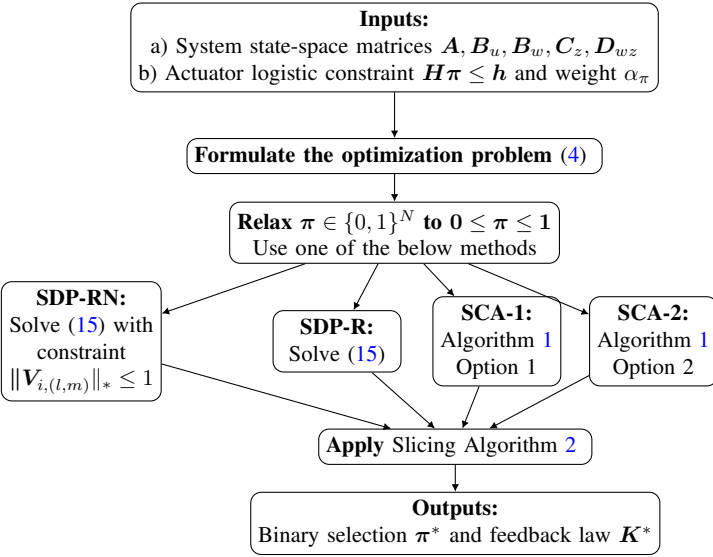
\begin{figure}
%	\centering
\hspace{-0.4cm}\begin{tikzpicture}[ node distance=1.2cm,
every node/.style={fill=white}, align=center, scale=0.85,every node/.style={scale=0.85}]
% Specification of nodes (position, etc.)
\node (start)             [activityStarts]              {\textbf{Inputs:} \\
	a) System state-space matrices
	$\mA, \mB_u, \mB_w, \m C_z, \m D_{wz}$\\
	b) Actuator logistic constraint $\mH\m\pi \leq \m h$ and weight $\alpha_{\pi}$
};
\node (firstproblem)     [process, below of=start, yshift=-0.5cm]          {\textbf{Formulate the optimization problem~\eqref{eq:StabUI2-all}}
};
%\node (firstproblemparam)     [process, left of=firstproblem, xshift=-3.5cm]          {\textbf{Set} $\m\alpha_{\pi}$
%};
\node (relax)  [process, below of=firstproblem]          {\textbf{Relax} $\m\pi\in\lbrace 0,1\rbrace^{N}$ \textbf{to} $\m 0 \leq \m\pi  \leq \m 1$ \\ Use one of the below methods};
\node (SDPR)     [process, below of=relax, yshift=-0.5cm, xshift = -0.9cm]          {\textbf{SDP-R: }\\
	Solve~\eqref{eq:StabUI-sdprelax-all}
};
\node (slicing)     [process, below of=SDPR, yshift=-0.45cm, xshift = 1.85cm]          {\textbf{Apply} Slicing Algorithm~\ref{algorithm:Slicing2}
};
\node (solution)     [process, below of=slicing, yshift=0cm]          {\textbf{Outputs:} \\
	Binary selection $\m\pi^*$ {and feedback law} $\m K^*$
};
\node (SDPRN)     [process, left of=SDPR, xshift=-2.8cm]          {\textbf{SDP-RN:} \\
	Solve~\eqref{eq:StabUI-sdprelax-all} with \\ constraint  \\ $\m \Vert \mV_{i,(l,m)}\Vert_* \leq 1$
};
\node (SCA1)     [process, right of=SDPR, xshift=1.3cm]          {\textbf{SCA-1:} \\
	Algorithm~\ref{algorithm:SCA}  \\ Option 1
};
\node (SCA2)     [process, right of=SCA1, xshift=1.25cm, yshift=0cm]          {\textbf{SCA-2:} \\
	Algorithm~\ref{algorithm:SCA}  \\ Option 2
};
\draw[->]             (start)  -- (firstproblem);
%\draw[->]             (firstproblemparam) -- (firstproblem);
\draw[->]             (firstproblem) -- (relax);
\draw[->]             (relax) -- (SDPR);
\draw[->]             (SDPR) -- (slicing);
\draw[->]             (slicing) -- (solution);
\draw[->]             (relax) -- (SDPRN);
\draw[->]             (SDPRN) -- (slicing);
\draw[->]             (relax) -- (SCA1);
\draw[->]             (relax) -- (SCA2);
\draw[->]             (SCA1) -- (slicing);
\draw[->]             (SCA2) -- (slicing);
\end{tikzpicture}	
\vspace{-0.12cm}
\caption{{Flow chart showing  the actuator selection and feedback control approach for the developed methods.}}
\label{fig:flowchart}
\vspace{-0.38cm}
\end{figure}

\normalcolor
\section{Numerical Tests}	\label{sec:numtests}
In this work, we develop different computational methods to solve the actuator selection problem with a focus on the $\mathcal{L}_{\infty}$ control metric~\eqref{eq:StabUI2-all}. We also showed that other control and estimation formulations can be formulated in the same fashion. The methods are summarized as follows.
%first four methods rely on relaxing the actuator integer constraints to box constraints.
\begin{itemize}
	\item \textit{SDP-R}: An SDP relaxation providing a lower bound to the optimal solution of the problem with BMIs; see~\eqref{eq:StabUI-sdprelax-all}.
	\item \textit{SDP-RN}: Same as SDP-R with the addition of the nuclear norm constraint to~\eqref{eq:StabUI-sdprelax-all}; see Corollary~\ref{corol:nucnorm}.
	\item \textit{SCA-1} and \textit{SCA-2}: Successive convex approximations producing upper bounds; see~\eqref{eq:StabUI-DC-all} and~\eqref{eq:SecondCA}.
	\item These four methods (SDP-R, SDP-RN, SCA-1, SCA-2) are based on relaxing the integer constraints, and then followed by a slicing algorithm that returns an integer actuator selection and an upper bound on the optimization problem with MIBMIs (Algorithm~\ref{algorithm:Slicing2}). {Fig.~\ref{fig:flowchart} shows a flowchart summarizing these four methods.}
	%	\end{itemize}
	\item \textit{Big-M}: The fifth method pertains to a formulation that transforms a problem with MIBMIs~\eqref{eq:StabUI2-all} into an MISDP via the Big-M method from Section~\ref{sec:bigM}.
\end{itemize}
All the simulations are performed using MATLAB R2016b running on 64-bit Windows 10 with Intel Core i7-6700 CPU with base frequency of 3.4GHz and 16 GB of RAM. YALMIP~\cite{lofberg2004yalmip} and its internal branch-and-bound solver are used as a modeling language and MOSEK~\cite{mosek2015mosek} is used as the SDP solver for all methods.
\begin{table*}[htbp]
	\vspace{-0.42cm}
	\caption{Final results after running Algorithm~\ref{algorithm:Slicing2} to recover the binary actuator selection and the actual system performance for the system with random network. The boldfaced numbers describe the method that outperformed other methods (the MISDP solver on YALMIP that implements the Big-M approach is terminated after 300 branch-and-bound iterations). For the Big-M method, the gap percentages are $1.2, 10.19, 25.31, 44.90,  47.33, 51.48,  52.63, 52.21,  53.54,  55.91$ for $N=5,10,\ldots, 50$.}
	\hspace{-1cm}
	\begin{tabular}{r|rrrrr|rrrrr|rrrrr}
		\cmidrule{2-16}    \multicolumn{1}{c}{\multirow{2}[4]{*}{$N$}} & \multicolumn{5}{c}{Performance Index $\sqrt{(\eta+1)\zeta}$} & \multicolumn{5}{c}{Total Activated Actuators $\sum_i^N \pi_i$} & \multicolumn{5}{c}{$f_{\mathrm{final}}=(\eta+1)\zeta+\sum_i^N \pi_i$} \\
		\cmidrule{2-16}    \multicolumn{1}{c}{} & \multicolumn{1}{c}{Big-M} & \multicolumn{1}{c}{SDP-RN} & \multicolumn{1}{c}{SDP-R} &\multicolumn{1}{c}{SCA-1} & \multicolumn{1}{c}{SCA-2} & \multicolumn{1}{c}{Big-M} & \multicolumn{1}{c}{SDP-RN} &  \multicolumn{1}{c}{SDP-R} & \multicolumn{1}{c}{SCA-1} & \multicolumn{1}{c}{SCA-2} & \multicolumn{1}{c}{Big-M} &\multicolumn{1}{c}{SDP-RN} &  \multicolumn{1}{c}{SDP-R } & \multicolumn{1}{c}{SCA-1} & \multicolumn{1}{c}{SCA-2} \\
		\midrule
		5     & 1.200 & 1.133 & 1.165 & 1.003 & 0.986 & 2     & 3     & 3     & 3     & 3     & \textbf{3.441} & 4.284 & 4.357 & 4.005 & 3.971 \\
		10    & 1.496 & 1.291 & 1.280 & 1.117 & 1.360 & 3     & 5     & 5     & 7     & 5     & \textbf{5.239} & 6.668 & 6.639 & 8.248 & 6.849 \\
		15    & 1.459 & 1.274 & 1.397 & 1.245 & 1.378 & 6     & 12    & 8     & 13    & 10    & \textbf{8.129} & 13.622 & 9.952 & 14.549 & 11.898 \\
		20    & 1.079 & 1.198 & 1.379 & 1.243 & 1.320 & 15    & 19    & 14    & 17    & 18    & 16.165 & 20.435 & \textbf{15.903} & 18.545 & 19.742 \\
		25    & 1.232 & 0.001 & 1.082 & 0.911 & 5.015 & 19    & 25    & 23    & 24    & 17    & \textbf{20.517} & 25.000 & 24.170 & 24.831 & 42.148 \\
		30    & 1.173 & 1.647 & 0.995 & 1.554 & 2.152 & 24    & 24    & 28    & 28    & 16    & 25.376 & 26.711 & 28.991 & 30.415 & \textbf{20.631} \\
		35    & 1.343 & 1.578 & 2.028 & 1.277 & 1.255 & 30    & 34    & 23    & 31    & 32    & 31.804 & 36.489 & \textbf{27.112} & 32.632 & 33.575 \\
		40    & 1.201 & 1.280 & 1.605 & 1.284 & 1.287 & 35    & 38    & 28    & 38    & 33    & 36.442 & 39.639 & \textbf{30.576} & 39.649 & 34.656 \\
		45    & 1.258 & 1.640 & 1.086 & 1.362 & 1.548 & 40    & 36    & 44    & 40    & 36    & 41.583 & {38.689} & 45.180 & 41.854 & \textbf{38.396} \\
		50    & 0.980 & 2.236 & 1.283 & 1.389 & 2.426 & 45    & 32    & 43    & 47    & 39    & 45.961 & \textbf{37.001} & 44.646 & 48.930 & 44.885 \\
	\end{tabular}%
	\label{tab:result4}%
	\vspace{-0.24cm}
\end{table*}%

%\vspace{-0.5cm}

\subsection{Simulated Dynamic Systems, Parameters, and Setup}
We use a randomly generated dynamic network from~\cite{MihailoSite,Motee2008} as a benchmark to test the presented methods.  {The random dynamic network has the following structure
\begin{align*}
	\dot{\m x}_i &= \underbrace{-\bmat{1 & 1 \\ 1 & 2}}_{\m A_i} \m x_i+ \sum_{i\neq j} e^{-\alpha (i,j)} \m x_j+ \bmat {0 \\ 1} (\m u_i+\m w_i), 
	\end{align*}
 where the coupling between nodes $i$ and $j$ is determined by the Euclidean distance $\alpha (i,j)$. These distances are unique for every $N$ and randomly generated inside a box of size $N/5$. Note that in these tests, we made the individual $\m A_i$ matrix for each subsystem to be stable (in comparison with~\cite{MihailoSite,Motee2008} where $\m A_i$ is unstable), so that the total number of unstable eigenvalues is smaller for the dynamic network ($\m A$ still has unstable eigenvalues). Keeping the same structure for the $\m A$ matrix as in~\cite{MihailoSite,Motee2008} yields the trivial solution of activating all actuators which is needed to guarantee an $\mathcal{L}_{\infty}$-stable performance---and hence the modification.}
%\subsection{}
We also use the following specific parameters and constraints.
%\begin{itemize} 
%	\item 
	The constraint $\m H \m\pi \leq \m h$ is represented as $\sum_{i=1}^{N}\m \pi_i \geq \lfloor N/4\rfloor$, where $\lfloor\cdot\rfloor$ denotes the floor function. We also set $\m\alpha_{\pi}^{\top} = \bmat{1,\ldots,1}$, that is all actuators have equal weight; $\alpha=1$ and $\eta=1$ (these constants appear in the LMIs).
%	\item
	 For SCA-1 and SCA-2, to obtain $\m S_0,\, \zeta_0$, and $\m Z_0$, we initialize by assuming that $\m \Pi_0 = 0.1 \m I_{n_{u}}$, and subsequently solving the $\mathcal{L}_{\infty}$ SDP with $\mS_0 \succeq \epsilon_1\mI_{n_{x}}$ and $\zeta_0 \geq \epsilon_1$, where $\epsilon_1 = 10^{-8}$. 
%	 
%	 \vspace{-0.2cm}
\subsection{Results and Comparisons}
\vspace{-0.2cm}
Table~\ref{tab:result4} depicts the results after applying Algorithm~\ref{algorithm:Slicing2} for SDP-R~\eqref{eq:StabUI-sdprelax-all}, SDP-RN, SCA-1~\eqref{eq:StabUI-DC-all}, and SCA-2~\eqref{eq:SecondCA}.  Algorithm~\ref{algorithm:Slicing2} is not applied to the Big-M solutions, as these solutions are originally binary.  Table~\ref{tab:result4} presents the performance index $\sqrt{(\eta+1)\zeta}$, the total activated actuators $\sum_i^N \pi_i$, and the objective function value $f_{\mathrm{final}}=(\eta+1)\zeta+\sum_i^N \pi_i$. The presented results for the Big-M method are for 300 iterations for the branch-and-bound solver of YALMIP. The maximum number of iterations is reached while the gap percentage is still between 1\% for  $N=5$ all the way to 56\% for $N=50$ (the gap, provided in the caption of Table~\ref{tab:result4}, increases as $N$ increases). Unfortunately, solving MISDPs would require weeks before the optimal solution (for larger values of $N$) is obtained and hence the choice of the default maximum iterations number of 300.   

The boldfaced numbers in $f_{\mathrm{final}}$ column in Table~\ref{tab:result4} depict the method with the smallest objective function value. The Big-M/MISDP formulation has been proposed before for SaA selection in linear systems~\cite{Chanekar2017,Taylor2017}. While Big-M yields the smallest $f_{\mathrm{final}}$ in some cases, the other methods (SDP-R, SDP-RN, SCA-1, SCA-2) yield better objective values, while requiring significantly less computational time---often orders of magnitude smaller than Big-M. In particular, Table~\ref{tab:ctun} shows the computational time (in seconds) for the five methods.
Since SDP-R solves only a single SDP, it is expected to be computationally more efficient than the other methods---this can be observed from Table~\ref{tab:ctun}. In addition, and since SCA-1 includes a smaller number of constraints and variables than SCA-2 (see Section~\ref{subsec:compareSCA}), the former requires less computational time in several of the simulations. However, there are instances where the SCAs require less computational time the than semidefinite relaxations (SDP-RN and SDP-R). The unifying theme here is that relaxing the integer constraints and using the convex approximations and relaxations is a good alternative to computationally costly MISDPs. {In addition, we emphasize that although some methods can yield the same number of activated actuators, the specific activated actuators from each method can be significantly different.} 
%More discussions and numerical experiments are included in~\cite{Taha2017d}.
%showing the performance of greedy algorithms, and the tests on the network of mass-spring systems. 
%
%\begin{figure}[h]
%	\vspace{-0.4cm}	
%	\centering
%	\includegraphics[scale=0.5]{ctun}
%	\caption{CPU time for the different methods with various values for the number of nodes $N$. These results correspond to the system of random network. }
%	\label{fig:comptime2}
%	\vspace{-0.4cm}
%\end{figure}

\begin{table}[t]
\vspace{-0.3cm}	\centering
	\caption{CPU time for the different methods with various values for the number of nodes $N$ for the random dynamic network.}
	\begin{tabular}{c|ccccc}
		$N$ & Big-M & SDP-RN & SDP-R & SCA-1 & SCA-2 \\
		\midrule
		5     & 3.92  & 1.84  & \textbf{1.45} & 2.10  & 1.84 \\
		10    & 87.47 & 3.73  & \textbf{1.36} & 2.97  & 2.58 \\
		15    & 369.49 & 14.18 & \textbf{3.36} & 10.49 & 8.86 \\
		20    & 1337.97 & 50.26 & \textbf{19.35} & 33.82 & 45.17 \\
		25    & 3774.93 & 142.73 & 90.88 & 120.51 & \textbf{80.12} \\
		30    & 9222.35 & 317.55 & \textbf{281.83} & 314.63 & 127.21 \\
		35    & 19760.87 & 853.73 & \textbf{303.23} & 615.68 & 674.92 \\
		40    & 41038.02 & 1901.40 & \textbf{822.92} & 1673.95 & 1258.57 \\
		45    & 76166.24 & 3103.24 & 3201.57 & 2695.33 & \textbf{2192.87} \\
		50    & 131035.62 & 4107.22 & 4441.03 & 5785.46 & \textbf{4096.90} \\
	\end{tabular}%
	\label{tab:ctun}%	\vspace{-0.4cm}
	\vspace{-0.4cm}
\end{table}%

%\vspace{-0.4cm}

\subsection{Extensions to Sensor Selection for Nonlinear Systems} In this paper, we only use the $\mathcal{L}_{\infty}$ control problem with actuator selection to exemplify how the proposed methods can provide insights into the solution of MIBMIs. We emphasize that all other CPS dynamics and control/estimation formulations (see Appendix~\ref{app:log}) with SaA selection can be solved using the methods we develop here.
For example, consider the sensor selection alongside the state estimator design problem for nonlinear systems $\dot{\m x}=\m A\m x+ \m B_u \m u+ \m \phi(\m x),\;\; \m y= \m \Gamma \m C \m x$ where $\m\phi(\m x)$ is the vector of nonlinearities with Lipschitz constant $\beta>0$ and $\m \Gamma$ is the binary sensor selection variable (cf.~\eqref{equ:CPSModelSaA}). By considering the last SDP in Table~\ref{tab:ObserverDesign}, the weighted sensor selection problem becomes:
\begin{align*}  
\min_{\m \Gamma, \m P, \m Y, \kappa} & \;\;\alpha_{\gamma}^{\top} \m \gamma \\
\st &\;\;\begin{bmatrix}
\m A ^{\top} \bm P + \bm  P \m A -\bm  Y \m \Gamma \m C \\-  \m C ^{\top} \m \Gamma \bm Y ^{\top} + \alpha \bm P +\kappa\beta^2 \m I & \bm P \\
\bm  P & - \kappa \m I\end{bmatrix}\preceq \m O\\ 
			&\;\; \m H \m \gamma \leq \m h, \; \m \gamma \in  \{0,1\}^N,
\end{align*}
which can be solved using the developed methods in the paper. This formulation yields observer gain $\m L^*=(\m P^{*})^{-1}\m Y^*$ that guarantees the asymptotic stability of the estimation error $\m e(t)=\m x(t)-\hat{\m x}(t)$ from $$\dot{\hat{\m x}}=\m A\hat{\m x}+\m B_u \m u + \m L^* (\m y - \hat{\m y}) + \m \phi(\hat{\m x}), \; \; \hat{\m y} =\m \Gamma^* \m C \hat{\m x}$$
with minimal number of sensors $\m \Gamma^*$. 

\normalcolor
%\vspace{-0.5cm}
\section{Summary and Future Work}
This paper puts forth a framework to solve SaA selection problems for uncertain CPSs with various control and estimation metrics. Given the widely popular SDP formulations of various control and estimation problems (without SaA selection), we present various techniques that aim to recover, approximate, or bound the optimal solution to the combinatorial SaA selection problem via convex programming. While the majority of prior art focuses on specific metrics or dynamics, the objective of this paper is to present a unifying framework that streamlines the problem of time-varying SaA selection in uncertain and potentially nonlinear CPSs. 
 
The developed methods in the paper have their limitation. First, the transition in the state-space matrices needs to be given before the time-varying actuator selection problem is solved. This narrows the scope of the actuator selection problem. In future work, we plan to study the actuator selection problem when the topological evolution is unknown, yet bounded. In particular, we plan to explore solutions to the actuator selection problem if the state-space matrix $\m A$ includes bounded perturbations that mimic the evolution in the CPS topology.

We plan to study the following related research problems in future work. (1) Simultaneous SaA selection for output feedback control problems: This problem produces more complex integer programs than MIBMIs. (2) Applications to selection of distributed generation in electric power networks with frequency-performance guarantees. (3) Customized branch-and-bound and cutting plane methods that can improve the performance of the Big-M method. {  (4) Theoretical analysis of the tightness of the lower and upper bounds resulting from the convex formulations in this paper for various CPSs.} 

%\vspace{-0.5cm}

\section*{Acknowledgments}
The authors would like to thank Johan L\"{o}fberg for his helpful comments and suggestions, and the anonymous reviewers and associate editor for their constructive criticism. 
%\balance
\bibliographystyle{IEEEtran}	\bibliography{bibfile}

\appendices 

\section{Actuator Selection: The Logistic Constraints}~\label{app:log}
The constraint $\m H\m \pi \leq \m h$ couples the selected actuators across time periods,
and is a linear logistic constraint that includes the following scenarios.
\begin{itemize} 
	\item Activation and deactivation of SaAs in a specific time-period $j$. For example, if actuator $i$ cannot be selected at period $j$, we set $\pi_i^j \leq 0$.
	\item If actuator $k$ is allowed to be selected only after actuator $i$ is selected at period $j$, we set
	$\pi_k^{j+1} \leq \pi_i^j,$ for $j=1,\ldots, T_f.$
	\item  If actuator $k$ must be deselected after actuator $i$ is selected at period $j$, we set
	$\pi_k^{j+1} \leq 1-\pi_i^j,$ for $j=1,...T_f$.
	\item Upper and lower bounds on the total number of active SaAs per period can be accounted for.
	\item Other constraints such as minimal number of required active actuators in a certain region of the CPS, and unit commitment constraints that are obtained from solutions day-ahead planning problems, can be included.
	%	\footnote{In power systems, for example, this could entail that a generator or distributed energy resource (i.e., an actuator) $k$ will become inactive for a certain periods of time. }
\end{itemize}
%\vspace{-0.4cm}
%
 \newpage
\section{Proofs of Various Results}~\label{app:allproofs}
%\vspace{-0.4cm}
%\section{Proof of Proposition~\ref{prop:CUppBd}}~\label{app:p1}
\begin{IEEEproof}[Proof of Proposition~\ref{prop:CUppBd}]
	To construct the upper bound~\eqref{eq:BMIuppbd}, the bilinear term is written as
	\begin{align}
		-\mB_u  \m\Pi  \mZ  -\mZ^{\top} \m\Pi  \mB_u^{\top} &= \frac{1}{2} \left[ \left( \mB_u  \m\Pi  - \mZ^{\top}\right) \left(\mB_u  \m\Pi  - \mZ^{\top}\right)^{\top}\right. \notag\\
		& \hspace{-0.6cm}- \left.  \left( \mB_u  \m\Pi  +\mZ^{\top}\right) \left(\mB_u  \m\Pi  + \mZ^{\top}\right)^{\top} \right]
		\label{eq:bilinDC}
	\end{align}
	The term $\left( \mB_u  \m\Pi  - \mZ^{\top}\right) \left(\mB_u  \m\Pi  - \mZ^{\top}\right)^{\top}$ is convex in $\mZ $ and $\m\Pi $ since it comes from an affine transformation of the domain of a convex function~\cite[Example~3.48]{boyd2004convex}. 
	The term 
	\[ \mathcal{H} (\m\Pi ,\mZ ):=- \left( \mB_u  \m\Pi  +\mZ^{\top}\right) \left(\mB_u  \m\Pi  + \mZ^{\top}\right)^{\top} \] 
	is concave in $\mZ $ and $\m\Pi $. 
	We can therefore invoke the fact that the first-order Taylor approximation of a concave function (at any point) is a global over-estimator of the function. Let $\m\Pi _0, \mZ_0 $ be the linearization point, and let $\mathcal{H}_{\mathrm{lin}} (\m\Pi ,\mZ ;\m\Pi _0, \mZ_0 )$ denote the linearization of $\mathcal{H} (\m\Pi ,\mZ )$ at the point $(\m\Pi _0, \mZ_0 )$. It holds that 
	\begin{equation}
	\mathcal{H} (\m\Pi ,\mZ ) \preceq \mathcal{H}_{\mathrm{lin}} (\m\Pi ,\mZ ;\m\Pi _0, \mZ_0 )
	\label{eq:LinUppBd}
	\end{equation}
	for all $\m\Pi _0, \mZ_0 $  and $\m\Pi , \mZ $.
	
	The linearization can be derived by substituting $\m\Pi  = \m\Pi_0  + (\m\Pi  -\m\Pi_0  )$ and $\mZ  = \mZ_0  + (\mZ  -\mZ_0  )$ into $\mathcal{H} (\m\Pi ,\mZ )$ and ignoring all second-order terms that involve  $ (\m\Pi  -\m\Pi_0  )$ and $ (\mZ  -\mZ_0  )$. The result is~\eqref{eq:HLinDef}.
	Combining~\eqref{eq:HLinDef} with~\eqref{eq:LinUppBd} and~\eqref{eq:bilinDC}, we conclude that the left-hand side of~\eqref{eq:StabUI2-MIBMI} is upperbounded as
	\begin{multline}
		\bmat{\mA \mS +\mS \mA^{\top}+\alpha\mS  \\-\mB_u  \m\Pi  \mZ  -\mZ^{\top} \m\Pi  \mB_u^{\top} &\mB_w  \\\mB_w^{\top} & -\alpha\eta \mI}  \\
		\preceq  \bmat{\mA \mS +\mS \mA^{\top}+\alpha\mS  \\ +\frac{1}{2}  \left( \mB_u  \m\Pi  - \mZ^{\top}\right) \left(\mB_u  \m\Pi  - \mZ^{\top}\right)^{\top} \\ +\frac{1}{2}  \mathcal{H}_{\mathrm{lin}} (\m\Pi ,\mZ ;\m\Pi _0, \mZ_0 ) &\mB_w  \\\mB_w^{\top} & -\alpha\eta \mI}.
		\label{eq:Cineq}
	\end{multline}
	This can be obtained using the fact that \[ \begin{bmatrix} \mA_1 & \mB \\ \mB^{\top} & \mC \end{bmatrix}\succeq \mO, \mA_2 \succeq \mA_1 \Longrightarrow  \begin{bmatrix} \mA_2 & \mB \\ \mB^{\top} & \mC\end{bmatrix} \succeq \m O \]
	which can be proved using the definition of positive semidefiniteness.
	Inequality~\eqref{eq:Cineq} holds for all $\m\Pi _0, \mZ_0 $ and $\m\Pi , \mZ $, and its left-hand side is $ \mathcal{C} (\m \Pi , \mZ  ;\m\Pi _0, \mZ_0 )$.
\end{IEEEproof}
%\section{Proof of Proposition~\ref{prop:DCconv}}~\label{app:p4}
\begin{IEEEproof}[Proof of Proposition~\ref{prop:DCconv}]
	Notice that problem~\eqref{eq:StabUI-DC2-all} has the same feasible set as~\eqref{eq:StabUI-DC-all}  (with $\m \Pi_0,\mZ_0$ replaced by $\m \Pi_{k-1},\mZ_{k-1}$). Corollary~\ref{corol:SCA1optval} establishes that its feasible set is a restriction of the one in~\eqref{eq:NLSDP-all}. It follows that $f(\m p_k)\geq L$, and $L_{k}^{(1)}\geq L$ holds because of the added regularizer in~\eqref{eq:StabUI-DC2-obj}.  The monotonicity of $\{f(\m p_k)\}_{k=1}^\infty$ follows from a corresponding result in~\cite[Lemma~4.2(c)]{dinh2012combining}. The sequence is thus monotone decreasing and bounded [the latter follows from the assumption on the boundedness of the feasible set of~\eqref{eq:NLSDP-all}]. It is a standard result in analysis that a bounded and monotone decreasing sequence has a limit. Therefore, $\hat{f}^{(1)}\geq L$ holds for the limit due to $f(\m p_k)\geq L$. The convergence result of part c) follows~\cite[Theorem~4.3]{dinh2012combining}. It is emphasized that the existence of at least one limit point is guaranteed by the boundedness of the feasible set. 
\end{IEEEproof}
%\section{Proof of Proposition~\ref{prop:K1UppBd}}~\label{app:p5}
\begin{IEEEproof}[Proof of Proposition~\ref{prop:K1UppBd}]
	Function $\mathcal{\m F}_{1}(\m p)$ is written as
	\begin{align*}
		\mathcal{\m F}_{1}(\m p) &= \mathcal{\mC}_0 + \mathcal{\mA}(\m p) + \mathcal{\mB}(\m p) \\ 
		&= \bmat { \mO & \mB_{w} \\ \mB_{w}^{\top} & - \alpha\eta \mI} + \bmat{\mA\mS + \mS\mA^{\top} + \alpha \mS & \mO \\ \mO & \mO} \\ & + \bmat {-\mB_u \m\Pi \mZ - \mZ^{\top}\m\Pi \mB_u^{\top} & \mO \\ \mO & \mO}.
	\end{align*}
	Substituting $ \m\Pi = \m\Pi_{k} + \Delta \m\Pi = \m\Pi_{k} + \m\Pi $ $\,- \,\m\Pi_{k} $ and $ \mZ = \mZ_{k}+ \Delta \mZ = \mZ_{k} +\mZ - \mZ_{k}$ into  $\mathcal{\mB}(\m z)$ yields
	\begin{align*}
		\mathcal{\mB}(\m p) &= \bmat{-\mB_{u}(\m\Pi_{k} + \Delta \m\Pi)(\mZ_{k} + \Delta \mZ) & \\ - (\mZ_{k} + \Delta \mZ)^{\top}(\m\Pi_{k} + \Delta \m\Pi)\mB_{u}^{\top} & \mO  \\ \mO & \mO },
	\end{align*}
	where $-\mB_{u}(\m\Pi_{k} + \Delta \m\Pi)(\mZ_{k} + \Delta \mZ)  
	= -\mB_{u}\m\Pi_{k}\mZ_{k} - \mB_{u}\m\Pi_{k} \Delta \mZ - \mB_{u}  \Delta \m\Pi \mZ_{k} -  \mB_{u} \Delta \m\Pi \Delta \mZ $ and 
	$-(\mZ_{k} + \Delta \mZ)^{\top}(\m\Pi_{k} + \Delta \m\Pi) \mB_{u}^{\top} = - \mZ_{k}^{\top}\m\Pi_{k}\mB_{u}^{\top}-\mZ_{k}^{\top}\Delta \m\Pi \mB_{u}^{\top} - \Delta \mZ ^{\top}\m\Pi_{k}\mB_{u}^{\top} - \Delta \mZ ^{\top}\Delta\m\Pi \mB_{u}^{\top}$.
	%\end{eqnarray*}
	
	Given this, $\mathcal{\mB}(\m p)$ can be rearranged as 
	\begin{align*}
		\mathcal{\mB}(\m p) &= \bmat{-\mB_{u}\m\Pi_{k}\mZ_{k}-\mZ_{k}^{\top}\m\Pi_{k}\mB_{u}^{\top}-\mB_{u}\m\Pi_{k}\Delta \mZ & \\ -\Delta \mZ ^{\top}\m\Pi_{k}\mB_{u}^{\top}-\mB_{u}\Delta \m\Pi \mZ_{k}-\mZ_{k}^{\top} \Delta \m\Pi \mB_{u}^{\top} & \mO  \\ \mO & \mO }  \\ 
		&\,\,\,\,\,\,\, +  \bmat {-\mB_{u}\Delta \m\Pi \Delta \mZ - \Delta \mZ ^{\top}\Delta \m\Pi \mB_{u}^{\top}& \mO  \\ \mO & \mO }.
	\end{align*}
	By combining and grouping these results, we obtain
	\begin{align*}
		\mathcal{\m F}_{1}(\m p)&=\bmat{-\mB_{u} \m\Pi_{k} \mZ_{k} - \mZ_{k} ^{\top}\m\Pi_{k}\mB_{u}^{\top} & \mB_{w} \\  \mB_{w} ^{\top} & -\alpha\eta \mI}\\ 	
		&+\bmat{\mA\mS+\mS\mA^{\top}+\alpha\mS -\mB_{u}\m\Pi_{k}\Delta \mZ & \\ - \Delta \mZ ^{\top}\m\Pi _{k}\mB_{u} ^{\top} - \mB_{u}\Delta\m\Pi \mZ_{k} - \mZ_{k} ^{\top}\Delta \m\Pi \mB_{u} ^{\top} & \mO \\ \mO & \mO} \\ 
		&+ \bmat{-\mB_{u}\Delta \m\Pi \Delta \mZ-\Delta \mZ ^{\top}\Delta \m\Pi \mB_{u} ^{\top} & \mO \\ \mO & \mO}.
	\end{align*}
	An upper bound for the last bilinear term for any $\boldsymbol{Q} \in \mathbb{S}^{n_u}_{++} $ is given as~\cite[Lemma 1]{LeeHu2016}
	\begin{align*}
		-\mB_{u}\Delta \m\Pi \Delta \mZ-\Delta \mZ ^{\top}\Delta \m\Pi \mB_{u} ^{\top} \preceq \\ 
		\mB_{u} \Delta \m\Pi \mQ \Delta \m\Pi \mB_{u}^{\top}+ \Delta \mZ ^{\top} \mQ^{-1}\Delta \mZ.
	\end{align*} 
	Combining the previous two results yields~\eqref{eq:K1UppBd}.
\end{IEEEproof}
%\section{Proofs of Lemmas~\ref{lem:MinQinvUppBd} and~\ref{lem:KsUppBd}}~\label{app:p6}
\begin{IEEEproof}[Proof of Lemma~\ref{lem:MinQinvUppBd}]
	Let $\mathcal{\mR}(\m x; \m x_k)$ be the first-order Taylor approximation of $-\mathcal{\mQ}(\m x)^{-1}$ computed around $\m x_k$. That is
	\begin{align}
		\label{eq:MinQinvUppBd-1}
		\mathcal{\mR}(\m x; \m x_k)= -\mathcal{\mQ}(\m x_k)^{-1}-[D\mathcal{\mQ}(\m x_k)^{-1}](\m x - \m x_k). 
	\end{align}
	By setting $\Delta\m x=\m x - \m x_k$, the differential $-[D\mathcal{\mQ}(\m x_k)^{-1}]\Delta\m x$ is given by~\cite{petersen2008matrix}
	\begin{align*}
		[D\mathcal{\mQ}(\m x_k)^{-1}]\Delta\m x &= -\mathcal{\mQ}(\m x_k)^{-1}[D\mathcal{\mQ}(\m x_k)]\Delta\m x\mathcal{\mQ}(\m x_k)^{-1} \\
		&= -\mathcal{\mQ}(\m x_k)^{-1}\sum_{i=1}^n \frac{\partial \mathcal{Q}(\m x_k)}{\partial x_{i}}\Delta x_i\mathcal{\mQ}(\m x_k)^{-1} \\
		%&\mspace{-100mu}= -\mathcal{\mQ}(\m x_k)^{-1} \mathcal{\mQ}(\m x)\mathcal{\mQ}(\m x_k)^{-1} +  \mathcal{\mQ}(\m x_k)^{-1} \mathcal{\mQ}(\m x_k) \mathcal{\mQ}(\m x_k)^{-1}.
		&= -\mathcal{\mQ}(\m x_k)^{-1} \mathcal{\mQ}(\m x)\mathcal{\mQ}(\m x_k)^{-1} +  \mathcal{\mQ}(\m x_k)^{-1}.
	\end{align*}
	Substituting the latter into~\eqref{eq:MinQinvUppBd-1} yields
	%-\mathcal{\mQ}(\m x)^{-1}|_{\m x_k} &= -\mathcal{\mQ}(\m x_k)^{-1} \\&\quad+\mathcal{\mQ}(\m x_k)^{-1}\sum_{i=1}^n \mathcal{\mQ}_i(\m x - \m x_k)\mathcal{\mQ}(\m x_k)^{-1} \\
	\begin{align*}
		\mathcal{\mR}(\m x; \m x_k) = -2\mathcal{\mQ}(\m x_k)^{-1}+\mathcal{\mQ}(\m x_k)^{-1}\mathcal{\mQ}(\m x)\mathcal{\mQ}(\m x_k)^{-1}.
	\end{align*}
	Since $\mathcal{\mQ}(\m x)$ is positive definite, then it follows that $-\mathcal{\mQ}(\m x)^{-1}$ is concave~\cite[Example~3.48]{boyd2004convex}. Because the first-order approximation of a concave function is a global over-estimator, we obtain~\eqref{eq:MinQinvUppBd}.
	%	\begin{align*}
	%		-\mathcal{\mQ}(\m x)^{-1} \preceq -2\mathcal{\mQ}(\m x_k)^{-1}+\mathcal{\mQ}(\m x_k)^{-1}\mathcal{\mQ}(\m x)\mathcal{\mQ}(\m x_k)^{-1}. 
	%	\end{align*} This verifies~\eqref{eq:MinQinvUppBd} and completes the proof.
\end{IEEEproof}
\begin{IEEEproof}[Proof of Lemma~\ref{lem:KsUppBd}]
	By linearizing $-\mQ^{-1}$ around a given $\mQ_k\in\mathbb{S}^{n_u}_{++}$, an upper bound on  $\mathcal{\m K}(\m p;\m p_{k},\m Q)$ can be derived as follows. Since $-\mQ^{-1}$ is concave in $\m Q$, then according to Lemma~\ref{lem:MinQinvUppBd}, the over approximation of $-\mQ^{-1}$ around $\mQ_{k}$ is $-2\mQ_{k}^{-1}+\mQ_{k}^{-1}\mQ \mQ_{k}^{-1}$.  Substituting this over approximation of $-\mQ^{-1}$ into $\mathcal{\m K}(\m p;\m p_{k},\m Q)$ and applying congruence transformation with $\mathrm{diag}(\mI,\mI,\mQ_{k},\mI)$ as the post and pre-multiplier yields~\eqref{eq:Ksdef}. The relation in~\eqref{eq:KsUppBd} is obtained due to the fact that $-\mQ^{-1} \preceq -2\mQ_{k}^{-1}+\mQ_{k}^{-1}\mQ \mQ_{k}^{-1}$.
\end{IEEEproof}

\begin{IEEEproof}[Proof of Proposition~\ref{prop:paramSCAconv}]
	The feasible set of problem~\eqref{eq:SecondCA} is a restriction of the one in~\eqref{eq:NLSDP-all} due to Proposition~\ref{prop:K1UppBd}, Lemma~\ref{lem:KLMI}, Lemma~\ref{lem:KsUppBd}. It therefore holds that  $f(\m p_k)\geq L$, and $L_{k}^{(2)}\geq L$ follows from  the addition of the regularizer in the objective. The monotonicity of $\{f(\m p_k)\}_{k=1}^\infty$ follows from a related result in~\cite[Lemma~6]{LeeHu2016}. The monotinicity and the boundedness imply the existence of the limit, similarly to Proposition~\ref{prop:DCconv}. The convergence in part c) is analogous to~\cite[Proposition~5]{LeeHu2016}. The existence of at least one limit point is ensured by the boundedness of the sequence $\{\m p_k\}_{k=1}^\infty$.
\end{IEEEproof}
\normalcolor

\section{Successive Convex Approximation Implementation}\label{app:scalgo}
Algorithm~\ref{algorithm:SCA} illustrates how the SCAs \eqref{eq:StabUI-DC2-all} and~\eqref{eq:SecondCA} can be solved sequentially until a maximum number of iterations ($\mathrm{MaxIter}$) or a stopping criterion defined by a tolerance ($\mathrm{tol}$) are met.
\begin{algorithm}
	\caption{Solving the successive convex approximations.}
	\label{algorithm:SCA}
	\begin{algorithmic}
		%		\STATE \textbf{initialize:} Set $k=0, \m \Pi_0 = \m I_{n_{u}}, \m Q_0 = \m I_{n_{x}}$; and $\zeta_0,\m S_0, \m Z_0$ to the solution of~\eqref{eq:StabUI-all}
		\STATE \textbf{input:} $\mathrm{MaxIterNum},\mathrm{tol}, k=0, \m \Pi_0 = \m I_{n_{u}}$
		%		 \STATE Initialization: $k=0$
		%		 , $J_{k+1}=\infty$
		\WHILE {$k < \mathrm{MaxIterNum} $}
		\STATE \texttt{Option 1:} Solve~\eqref{eq:StabUI-DC2-all}
		\STATE
		 \texttt{Option 2:} Solve~\eqref{eq:SecondCA}
		\IF{ $|\hat{L}_{k}^{(1)\,\mathrm{or}\,(2)}-\hat{L}_{k-1}^{(1)\,\mathrm{or}\,(2)}| < \mathrm{tol}$ }
		\STATE \textbf{break}
		\ELSE 
		\STATE $k\leftarrow k+1$
		\ENDIF		\ENDWHILE
		\STATE \textbf{output:} \{$\m S^\star, \zeta^\star, \m Z^\star, \m \Pi^\star \} \leftarrow$ $\{\m S^k, \zeta^k, \m Z^k, \m \Pi^k \}$
	\end{algorithmic}
\end{algorithm}
\vspace{-0.4cm}
 \section{Recovering the Binary Selection}~\label{sec:slicing}
The solutions obtained from~\eqref{eq:StabUI-sdprelax-all}, \eqref{eq:StabUI-DC2-all}, and~\eqref{eq:SecondCA} produce a noninteger solution for the actuator selection problem. Since the objective is to determine a binary selection for the actuators, we present in this section a simple \textit{slicing routine} that returns a binary selection given the solutions to the optimization problems in Sections~\ref{sec:relaxMIBMIs} and \ref{sec:succapproxMIBMIs}.

The slicing routine is presented in Algorithm~\ref{algorithm:Slicing2}. Since the objective of the $\mathcal{L}_{\infty}$ problem is to find a feedback gain $\m K = \m Z \m S^{-1}$ that renders the closed-loop system stable, the slicing algorithm ensures that the spectrum $\Lambda(\m A_{\mathrm{cl}})$ of the closed-loop system matrix $\m A_{\mathrm{cl}} = \m A -\m B_u\m \Pi \m K$ lies on the left-half plane. 

The slicing routine takes as an input the real-valued solution to the actuator selection $\m \Pi^*$ with $\m \pi_{i}^* \in [0,1]$. First, the entries of $\m \pi^*$ are sorted in decreasing order, and the minimum $s$-actuator selection is obtained such that the logistic constraints $\m H \m \pi \leq \m h$ are satisfied, given that $\m \pi \in \{0,1\}^N$. This ensures that we start the slicing algorithm from the minimum number of actuators, while still satisfying all of the actuator-related constraints in~\eqref{eq:StabUI2-all}. The algorithm proceeds by activating the $s$-highest ranked actuators, followed by solving the $\mathcal{L}_{\infty}$ SDP~\eqref{eq:StabUI2-OF}--\eqref{eq:StabUI2-zetaS} for $\m Z$ and $\m S$. Then, the maximum real part of the eigenvalues of $\m A_{\mathrm{cl}}$, namely $\lambda_m$, is obtained. If $\lambda_m < 0$, the algorithm exits returning the actuator selection $\m \Pi_s$ and the associated feedback gain. 

The algorithm allows the addition of other user-defined requirements, such as a minimum performance index $\zeta$ or a maximum $\lambda_m$, which can guarantee a minimal distance to instability. It  can also be generalized to other control or estimation problems.  {Notice that  Algorithm~\ref{algorithm:Slicing2} terminates when $\lambda_m < 0$ and the SDP~\eqref{eq:StabUI2-OF}--\eqref{eq:StabUI2-zetaS} is solved. These conditions ensure  by definition that the system is controllable for the resulting binary actuator combination. In short, the slicing algorithm guarantees the controllability of the system.  }

%\vspace{-0.251cm}
%\vspace{-0.251cm}

The actuator selection and associated control law returned by Algorithm~\ref{algorithm:Slicing2} yield an upper bound $U$ to the optimal value of the actuator selection problem~\eqref{eq:StabUI-all}. 
%\begin{prop}
%	\label{prop:slicing}
%	Let $\m p_s$ denote the solution obtained from Algorithm~\ref{algorithm:Slicing2}. It holds that $f(\m p_s) = U \geq f^*$.
%\end{prop}
%%\begin{IEEEproof}[Proof of Proposition~\ref{prop:slicing}]
%%By construction of Algorithm~\ref{algorithm:Slicing2}, the point $\m p_s$ is feasible for~\eqref{eq:StabUI-all} since the inequality constraints on the actuator selection are satisfied. Thus, it follows that $f(\m p_s)=U \geq f^*$.
%\end{IEEEproof}

 \begin{algorithm}
 	\caption{A Slicing Algorithm to Recover the Integer Selection from \eqref{eq:StabUI-sdprelax-all}, \eqref{eq:StabUI-DC2-all}, and~\eqref{eq:SecondCA}}\label{algorithm:Slicing2}
 	\begin{algorithmic}
 		%		\STATE \textbf{initialize:} $F \leftarrow \emptyset, k = 0$
 		\STATE \textbf{input:} $\m \Pi^*$ from Algorithm~\ref{algorithm:SCA}, set $ \lambda_{m}= \infty$
 		\STATE Sort $\m \pi^*$ in a decreasing order
 		% 		\STATE  as the least number of actuators that satisfy $\m H \m\pi \leq \m h$
 		\STATE $s = \mathop{\mathrm{minimum}}\limits_{\m \pi  \in   \{0,1\}^N, \m H \m\pi \leq \m h } \quad \m 1_{N}^{\top}\m \pi $
 		\WHILE {$ \lambda_{m}\geq 0$}
 		\STATE Activate the  $s$-highest ranked actuators in $\m \pi$
 		\STATE Obtain $\m {\Pi}_s = \mathrm{blkdiag}(\pi_1\m {I}_{n_{u_{1}}},\ldots,\pi_N\m {I}_{n_{u_{N}}})$
 		\STATE Given $\m \Pi=\m\Pi_s$, solve the SDP~\eqref{eq:StabUI2-OF}--\eqref{eq:StabUI2-zetaS} for $\m Z$ and $\m S$
 		\STATE $ \lambda_{m}= \max(\mathrm{real}(\lambda))$ where $\lambda \in \Lambda(\m {A}-\m {B_u}\m {\Pi}_s\m Z\m S^{-1})$
 		\STATE $s \leftarrow s+1$
 		\ENDWHILE
 		\STATE \textbf{output:} $\m \Pi_s^*, \m K^* = \m Z^* (\m S^*)^{-1}$
 	\end{algorithmic}
 \end{algorithm} 
  \onecolumn
%\vspace*{-0.7cm}

\section{Various Controller and Observer Designs via SDP Formulations for Different Metrics and Dynamics}~\label{app:cont}
\vspace{-0.4cm}
\begin{table}[h]
	\vspace{-0.35cm}	
	\centering 
	\caption{Controller Design for various CPS dynamics and objectives via SDP formulations~\cite{zhou1996robust,boyd1994linear,pancake2000analysis}.}
%	\vspace{-0.2cm}
	\label{tab:Controllers}
	\renewcommand{\arraystretch}{1.5}
	\begin{tabular}{ | l | l | l }
		\hline %inserts horizontal line
		\textbf{	CPS Dynamics, Metrics, \& Design Objective }&\textbf{Control Design via SDPs }\\[0.5ex] 
		\hline 
		\hline
		$\begin{array}{l}
		\textit{Stabilization of Linear Systems} \\
		\dot{\m x}  = \mA\m x +\mB_u \m u ,\; \text{Variable:}\,\, $ $\mS$ $
		\end{array}$ & $\begin{array} {lll} \mathrm{find} & \m S \;\;\;
		\st &		\m A\mS+\mS \m A^{\top} \preceq \m B_u \m B_u^{\top} \end{array}$ \\ 
		\hline 
		$\begin{array}{lll}
		\textit{Robust}\; \mathcal{L}_{\infty}\;\textit{Control of Uncertain Linear Systems} \\
		\dot{\m x} = \mA\m x +\mB_u \m u +\mB_w \m w  \\
		\m z  = \mC_z \m x  +\mD_{vw} \m w  \\
		\m u  = - \mK \m x  =-\mZ\mS^{-1} \m x  ,\;
		\text{Variables:}\,\, \mZ,\mS,\zeta 
		\end{array}$ & 
		$\begin{array} {l} 	\min\,\,\, (\eta+1)\zeta \\
		\st  \bmat{\mA\mS+\mS\mA^{\top}-\mB_u \mZ -\mZ^{\top} \mB_u^{\top} +\alpha\mS &\mB_w \\\mB_w^{\top} & -\alpha\eta \mI} \preceq \m {O} \\
		\hspace{0.45cm}\bmat{-\mS & \mO & \mS\mC_z^{\top}\\
			\mO & -\mI & \mD_{wz}^{\top}\\
			\mC_z\mS& \mD_{wz} &-\zeta\mI} \preceq \m {O},  \zeta > 0 \end{array}$  \\ 
%		\hline 
%		$\begin{array}{lll}
%		\textit{Stabilization through Output Feedback Control } \\
%		\dot{\m x} = \mA\m x +\mB_u \m u  \\
%		\m y  = \mC \m x  \\
%		\m u  = - \mK \m y  =\mM^{-1}\mN \m y  ,\;
%		\text{Variables:}\,\, \mM,\mN,\mS
%		\end{array}$ & $\begin{array} {lll} \mathrm{find} & \m S, \m N, \m M \\
%		\st &  \mA^{\top}\mS+\mS\mA-\mC^{\top}\mN^{\top}\mB_u^{\top}-\mB_u\mN\mC \preceq \mO ,\;\;\mB_u\mM=\mS\mB_u\end{array}$  \\ 
		\hline 
		$\begin{array}{lll}
		\textit{LQR Control: Minimizing State and Input Energy} \\
		\min  \,\,\, \,\, \mathbb{E} \int_{t_0}^{\infty} \m x^{\top}(\tau)\mQ\m x(\tau)+\m u^{\top}(\tau)\mR \m u(\tau) d\tau\\
		\st \,\,\, \,\, \dot{\m x}  = \mA\m x +\mB_u \m u + \m w, \m w \sim \mathcal{N}(\m 0, \m W)
		\\
		\m u  = -\mR^{-1}\mB_u^{\top}\mS^{-1}\m x  ,\;
		\text{Variables:}\,\, \mY,\mS 
		\end{array}$ & $\begin{array} {lll} 
		\min  & \,\,\,\,\mathbf{trace}(\m W\m S^{-1})\\
		\st & \bmat{\mA\mS+\mS\mA^{\top}+\mB_u \mY + \mY^{\top} \mB_u^{\top} &\mS & \mY \\\mS & -\mQ^{-1} & 0 \\ \mY^{\top} & 0 & -\mR^{-1}} \preceq \m {O} \end{array}$  \\ 
		\hline
		$\begin{array}{lll}
		\textit{Stabilization of Time-Delay Systems} \\
		\dot{\m x} = \mA\m x +\mB_u \m u  + \sum_{i=1}^{L} \m A_i \m x(t-\tau_i)\\
		%			\m y  = \mC \m x 
		%			\\
		\m u =\mK\m x =\mZ\mS^{-1}\m x ,\;
		\text{Variables:}\,\, \mZ,\mS, \m S_1, \ldots, \m S_L
		\end{array}$ & $\begin{array} {l} 		\mathrm{find} \,\, \mZ,\mS, \m S_1, \ldots, \m S_L \\
		\st  \bmat{\mA\mS+\mS\mA^{\top}+\mB_u \mZ +\mZ^{\top} \mB_u^{\top} + \sum_{i=1}^{L} \m S_i &\m A_1 \m S & \ldots & \m A_L \m S\\ 
			\m S \m A_1^{\top}  & -\m S_1 & \ldots & \m O\\
			\vdots & \vdots & \ddots & \vdots \\
			\m S \m A_L^{\top} & \m O & \ldots & -\m S_L} \prec \m {O} \\
		\hspace{0.5cm} \end{array}$  \\
		\hline 
	\end{tabular}
%	\vspace{-0.7cm}
\end{table}

\begin{table}[!h]
	\centering 
	\caption{Estimators Design for CPSs via  SDP formulations; see references \cite{scherer2000linear,chakrabartystate,rajamani1998observers}. Other similar formulations are omitted for brevity.}
%	\vspace{-0.22cm}
	\label{tab:ObserverDesign}
	\renewcommand{\arraystretch}{1.5}
	\begin{tabular}{|l|l|l|}
		\hline %inserts horizontal line
		\textbf{ Dynamics \& Estimation Objective } & \textbf{Estimator Dynamics }& \textbf{Estimator Design via SDPs} \\[0.5ex] 
		\hline 
		\hline
		$\begin{array}{l}
		%		\textit{Nominal Linear System}\\
		\dot{\m x} = \m A  \m x + \m B_u   \m u \\
		\m y  =  \m C  \m x + \m D_u  \m u  \\ \text{Variables:}\,\,  \bm P, \bm Y
		%			\,\,\text{Reference \cite{luenberger1966observers}.}
		\end{array}$ & $\begin{array} {lcl} 
		\dot{\hat{  \m x}} =  \m A\hat{ \m  x} +\m  B_u   \m u + \bm L ( \m  y -\hat{\m   y} )\\
		\hat{\m   y}  = \m C \hat{ \m  x}  + \m D_u   \m u \\ \bm L= \bm P^{-1}\bm  Y  \end{array}$ & $\begin{array} {lcl}  \mathrm{find} &  \bm Y \\
		\st
		&\m  A^{\top} \bm P+ \bm P  \m A - \m C^{\top} \bm Y^{\top}- \bm Y \m C +\alpha \bm P \preceq \m O\\ 
		%			& \bm P \succ \m O
		\end{array}$\\ 
		\hline			$\begin{array}{l}
		%			\textit{Linear System with UIs}\\
		\dot{\m x} = \m A  \m x + \m B_u   \m u + \m B_w  \m w \\
		\m y  =  \m C  \m x + \m D_u  \m u  \\ \text{Variables:}\,\,  \bm P,\bm Y,\mu
		%			\,\,\text{Reference~\cite{scherer2000linear}.}
		\end{array}$ & $\begin{array} {lcl}  	\dot{\hat{\m x}} =\m  A \hat{\m x} + \m B_u \m u + \bm L(\m y -\hat{\m y} ) \\
		\hat{\m y}  = \m C \hat{\m x}  + \m D_u \m u  \\
		\bm L=\bm P^{-1}\bm Y \end{array}$ & 
		$\begin{array}{lcl}
		\min &  \mu >0 \\
		\st	&  \begin{bmatrix}
		\m A^\top \bm P + \bm P\m A - \\\m C^\top \bm Y^\top -\bm Y\m C-\mu  \m C^\top \m C & \bm P \m B_w  \\ 
		\m B_w^{\top} \bm P & - \mu  \m I 
		\end{bmatrix} \preceq \m O\\
		&  \;  \mu  > 0
		\end{array} $\\
		\hline
		$\begin{array}{l}
		\dot{\m x} = \m A  \m x + \m B_u   \m u + \m B_w   \m w \\
		\m y  =  \m C  \m x + \m D_u  \m u   +  \m D_v   \m v \\ \text{Variables:}\,\,  \bm P,\bm  Y,  \mu
		%			 \text{~Reference~\cite{chakrabartystate}}
		\end{array}$ & $\begin{array} {lll} 		\dot{\hat{ \m x}} =  \m A\hat{  \m x} + \m B_u   \m u +\bm  L (  \m y -\hat{\m   y} )\\
		\hat{ \m  y}  =\m  C   \hat{\m x}  + \m D_u   \m u \\ \bm L=\m  P^{-1}\bm  Y  \end{array}$ & $\begin{array} {lcl} 	 \min & \mu  \\\st& \begin{bmatrix}
		\m A ^{\top} \bm P +  \bm P \m A -  \bm Y \m C \\-  \m C ^{\top} \bm Y ^{\top} + \alpha  \bm P&  \bm P \m B_w&  \bm Y \m D_v\\
		\m B_w ^{\top} \bm P& -\alpha  \m I& \m O\\
		\m D_v ^{\top} \bm Y ^{\top}& \m O& -\alpha  \m I
		\end{bmatrix} \preceq \m O\\  & \bm P \succeq  \mu \m I, \;\; \mu >0 \end{array}$\\ 
	\hline
		$\begin{array}{l}
		%		\textit{Lipschitz Nonlinear System}\\
		\dot{\m x} = \m A  \m x + \m B_u   \m u +  \m  \phi(\m x) \\
		\m y  =  \m C  \m x + \m D_u  \m u  \\ \text{Variables:}\,\, \bm  P, \bm Y, \kappa
		%			\,\,\text{Reference \cite{rajamani1998observers}.}
		\\ \beta\,\, \text{is the Lipschitz constant \cite{rajamani1998observers}.} 			\end{array}$ & $\begin{array} {lcl} 	\dot{\hat{\m   x}} = \m   A \hat{  \m x} + \m B_u   \m u +  {{\m \phi}}(\hat{\m   x}) \\-  \bm L(  \m y  - \hat{  \m y} )\\
		\hat{  \m y}  = \m C \hat{ \m  x}+ \m D_u   \m u\\
		\bm L= \bm P^{-1}\bm  Y  \end{array}$ & $\begin{array} {lcl}  \mathrm{find} &  \bm Y,\kappa >0 \\
		\st&\begin{bmatrix}
		\m A ^{\top} \bm P + \bm  P \m A -\bm  Y \m C \\-  \m C ^{\top} \bm Y ^{\top} + \alpha \bm P +\kappa\beta^2 \m I & \bm P \\
		\bm  P & -\bm \kappa I\end{bmatrix}\preceq \m O\\ 
		%			& \; \kappa > 0
		\end{array}$\\ 
		\hline
	\end{tabular}
\end{table}

\end{document}